\documentclass[12pt]{article}
\usepackage{amsfonts}

\begin{document}

\def\be{\begin{equation}}
\def\ee{\end{equation}}
\def\bea{\begin{eqnarray}}
\def\eea{\end{eqnarray}}
\def\non{\nonumber}
%
%
\newtheorem{thm}{Theorem}[section]
\newtheorem{lem}[thm]{Lemma}
\newtheorem{prop}[thm]{Proposition}
\newtheorem{cor}[thm]{Corollary}
\newenvironment{dfn}{\medskip\refstepcounter{thm}
\noindent{\bf Definition \thesection.\arabic{thm}\ }}{\bigskip}
\newenvironment{example}{\medskip\refstepcounter{thm}
\noindent{\bf Example \thesection.\arabic{thm}\ }}{\medskip}
\newenvironment{proof}{\medskip\noindent{\it Proof.}}{\hfill
\rule{.5em}{.8em}\medskip}

\newenvironment{exmp}{\refstepcounter{exmp}
\medskip\noindent{\bf Example \arabic{exmp}.\ }}{\medskip}
\def\ms#1{\vert#1\vert^2}
\def\bms#1{\bigl\vert#1\bigr\vert^2}
\def\nm#1{\Vert #1 \Vert}
\def\md#1{\vert #1 \vert}
\def\bnm#1{\bigl\Vert #1 \bigr\Vert}
\def\bmd#1{\bigl\vert #1 \bigr\vert}
\def\Bbb{\mathbb}
\def\goth{\mathfrak}
\def\k{K\"ahler\ }
\def\kk{K\"ahler}
\def\hk{hyperk\"ahler\ }
\def\Hk{Hyperk\"ahler\ }
\def\hkk{hyperk\"ahler}
\def\Hkk{Hyperk\"ahler}
\def\qn{quaternion}
\def\qns{quaternions}
\def\qc{quaternionic\ }
\def\Qc{Quaternionic\ }
\def\qcc{quaternionic}
\def\Qcc{quaternionic}
\def\sc{symplectic\ }
\def\Sc{Syplectic\ }
\def\sy{symplectically }
\def\qk{\qc \k}
\def\Qk{\Qc \k}
\def\qkk{\qc \kk}
\def\Qkk{\Qc \kk}
\def\hc{hypercomplex\ }
\def\Hc{Hypercomplex\ }
\def\hcc{hypercomplex}
\def\Hcc{Hypercomplex}
\def\dim{\mathop{\rm dim}}
\def\mod{\mathop{\rm mod}}
\def\Re{\mathop{\rm Re}}
\def\Im{\mathop{\rm Im}}
\def\im{\mathop{\rm im}}
\def\Ker{\mathop{\rm Ker}}
\def\rk{\mathop{\rm rk}}
\def\Coker{\mathop{\rm Coker}}
\def\ind{\mathop{\rm ind}}
\def\End{\mathop{\rm End}}
\def\Aut{\mathop{\rm Aut}}
\def\Vect{\mathop{\rm Vect}}
\def\curl{\mathop{\rm curl}}
\def\Div{\mathop{\rm div}}
\def\Hom{\mathop{\rm Hom}}
\def\Hol{\mathop{\rm Hol}}
\def\Ric{\mathop{\rm Ric}}
\def\vol{\mathop{\rm vol}}
\def\id{\mathop{\rm id}}
\def\ah{{\rm A}$\Bbb H$}
\def\sh{{\rm S}$\Bbb H$}
\def\hh{$\Bbb H$}
\def\h{$\Bbb H$\ }
\def\ihh{$\Bbb I$}
\def\ih{$\Bbb I$\ }
\def\H{\ifmmode{\mathbin{\Bbb H}}\else $\Bbb H$\fi}
\def\R{\ifmmode{\mathbin{\Bbb R}}\else $\Bbb R$\fi}
\def\Z{\ifmmode{\mathbin{\Bbb Z}}\else $\Bbb Z$\fi}
\def\C{\ifmmode{\mathbin{\Bbb C}}\else $\Bbb C$\fi}
\def\CP{\ifmmode{\mathbin{\Bbb CP}}\else $\Bbb CP$\fi}
\def\HP{\ifmmode{\mathbin{\Bbb HP}}\else $\Bbb HP$\fi}
\def\I{\ifmmode{\mathbin{\Bbb I}}\else $\Bbb I$\fi}
\def\Q{\ifmmode{\mathcal{Q}}\else $\mathcal{Q}$\fi}
\def\O{\ifmmode{\mathcal{O}}\else $\mathcal{O}$\fi}
\def\a{{\mathbf a}}
\def\b{{\mathbf b}}
\def\x{{\mathbf x}}
\def\y{{\mathbf y}}
\def\p{\prime}
\def\d{\dagger}
\def\br{\buildrel}
\def\ov{\overline}
\def\ot{\otimes}
\def\ra{\rightarrow}
\def\lra{\longrightarrow}

\def\DD{\overline{{\mathcal{D}}}}
\def\D'{{\mathcal{D}}}
\def\Dol{{\overline{\partial}}}
\def\ovD{\overline{D}}
\def\bd{\bar{\delta}}

\def\longra{\longrightarrow}
\def\t{\times}
\def\ha{{\textstyle{1\over2}}}
\def\oth{{\mathchoice{\textstyle\otimes_{\Bbb H}}{\textstyle\otimes_
{\Bbb H}}{\scriptstyle\otimes_{\scriptscriptstyle\Bbb H}
}{\scriptscriptstyle\otimes_{\scriptscriptstyle\Bbb H}}}}
\def\bigoth{{\textstyle\bigotimes}_\H}
\def\Lambh{{\textstyle\Lambda}_\H}
\def\Sh{S_\H}
\def\Hah{{\textstyle\Hom_{{\mathrm A}\H}}}
\def\Aah{{\textstyle\Aut_{{\mathrm A}\H}}}
\def\op{\oplus}
\def\AA{{\cal A}}
\def\Sp{\ifmmode{\mathrm{Sp}}\else Sp\fi}
\def\GL{\ifmmode{\mathrm{GL}}\else GL\fi}
\def\SU{\ifmmode{\mathrm{SU}}\else SU\fi}
\def\Gg{{\goth g}}
\def\u{{\goth u}}
\def\su{\goth{su}}
\def\sp{\goth{sp}}
\def\sl{\goth{sl}}
\def\gl{\goth{gl}}
\def\Fu{F^\uparrow}
\def\Fd{F^\downarrow}


\title{A Dolbeault-type Double Complex on Quaternionic Manifolds}
\author{Dominic Widdows} \maketitle


\begin{abstract}

It has long been known that differential forms on a complex manifold 
$M^{2n}$ can be
decomposed under the action of the complex structure to give the Dolbeault
complex. This paper presents an analogous double complex for a quaternionic
manifold $M^{4n}$ using the fact that its cotangent space $T^*_mM$ is 
isomorphic to the quaternionic vector space $\H^n$. This defines an action of 
the group Sp(1) of unit quaternions on $T^*M$, which induces an action of Sp(1) 
on the space of $k$-forms $\Lambda^kT^*M$. A double complex is obtained by 
decomposing $\Lambda^kT^*M$ into irreducible representations of Sp(1), resulting 
in new `quaternionic Dolbeault' operators and cohomology groups. 

Links with previous work in quaternionic geometry, particularly the differential 
complex of Salamon and the q-holomorphic functions of Joyce, are demonstrated.

\end{abstract} 

\section*{Introduction} \label{Introsec}

This paper describes a new double complex of differential forms on hypercomplex
or quaternionic manifolds.  This is to date the clearest quaternionic version of
the more familiar Dolbeault complex, used throughout complex geometry.  It is
hoped that readers familiar with complex geometry will find the new ideas both
natural and familiar.  An important thread will be to examine the role played by
the group U(1) of unit complex numbers in complex geometry:  complex geometry is
then translated into quaternionic geometry by replacing U(1) with the group
Sp(1) of unit quaternions.  This approach is rewarding because of the simplicity
of the representations of the groups U(1) and Sp(1), which makes their action
much easier to understand than those of the groups $\GL(n,\C)$, $\GL(n,\H)$ and
$\Sp(1)\GL(n,\H)$.

The paper is arranged as follows.  Section \ref{BackSec} contains background
material from \mbox{complex} and quaternionic geometry, reviewing the definition of a
complex manifold and considering quaternionic analogues.  At least two
definitions are currently recognised.  Following the work of Salamon
\cite{SalQgeom}, the term `quaternionic manifold' is used for a 
\mbox{manifold}
possessing a torsion-free $\Sp(1)\GL(n,\H)$-structure.  The more restricted class
of `hypercomplex manifolds' refers to those with a torsion-free 
$\GL(n,\H)$-structure.  \mbox{Hypercomplex} manifolds possess global anticommuting
complex structures (since \linebreak 
$\GL(n,\H)\subset\GL(2n,\C))$.  Quaternionic manifolds
possess a 2-sphere bundle of local anticommuting almost-complex structures, but
whilst this bundle is invariant it may have no global sections, so there may be
no global complex structures.  The Riemannian versions of these structures
(K\"{a}hler, \hk and \qc K\"{a}hler manifolds) are briefly described.  This
section also recalls the most basic facts about Sp(1)-representations, their
weights and tensor products, which we later use to decompose forms on
quaternionic manifolds. Associated bundles are described. 

Section \ref{Dolsec} reviews the decomposition of differential forms in complex
geometry, and the resulting Dolbeault complex.  (This material will be familiar
to most readers.)  The space $\Lambda^{p,q}M$ of $(p,q)$-forms on a complex
manifold $(M,I)$ is a representation of the Lie algebra $\u(1)$ via the induced
action of the complex structure $I$ on $\Lambda^kT^*M$, a point of view which
adapts well to quaternionic geometry.  As well as the standard decomposition of
complex-valued forms, we also describe the less familiar decomposition of
real-valued forms.  The resulting `real Dolbeault complex' is even more closely
akin to the new quaternionic complex because they both form isosceles triangles
(as opposed to the diamond configuration of the standard Dolbeault complex).

Section \ref{Backgroundsec} surveys previous work on the decomposition of
differential forms in quaternionic geometry.  Kraines \cite{YK} and Bonan
\cite{Bonan} obtained a decomposition of forms on \qc K\"{a}hler manifolds by
taking successive exterior products with the fundamental 4-form.  Swann
\cite{Swann} considered decomposition of these forms as
$\Sp(1)\Sp(n)$-representations.  In the non-Riemannian setting, Salamon
\cite{SalQgeom} used the coarser decomposition of forms on quaternionic
manifolds into $\Sp(1)\GL(n,\H)$-representations, resulting in a differential
complex which forms the top row of the new double complex.  Much of this algebra
and geometry can be inferred from Fujiki's comprehensive article \cite{Fujiki},
which describes much of the theory underlying this whole area of research.

The heart of this paper is Section \ref{DCsec} which constructs the analogue of
the Dolbeault complex in quaternionic geometry.  The complex structure $I$ is
replaced by the (possibly local) almost complex structures $I$, $J$ and $K$, the
group $\mathrm{U}(1)$ of unit complex numbers is replaced by the group $\Sp(1)$
of unit quaternions, and the Lie algebra $\u(1)=\langle I\rangle$ is replaced by
$\sp(1)=\langle I,J,K\rangle$.  Despite the possible lack of global complex
structures on quaternionic manifolds, the Casimir operator $I^2+J^2+K^2$ still
makes invariant sense.  It follows that the decomposition of $\Lambda^kT^*M$
into irreducible $\Sp(1)$-representations is also invariant.  A straightforward
calculation using weights leads to the following result (Proposition
\ref{LambdaMdecomp}):

\[\Lambda^kT^*M\cong\bigoplus_{r=0}^k
       \left[{2n\choose\frac{k+r}{2}}{2n\choose\frac{k-r}{2}}-
       {2n\choose\frac{k+r+2}{2}}{2n\choose\frac{k-r-2}{2}}\right]V_r, \] where
       $M^{4n}$ is a quaternionic manifold, $V_r$ is the Sp(1)-representation
       with highest weight $r$, and $r\equiv k \bmod 2$.  It follows from the
       Clebsch-Gordon formula $V_r\ot V_1\cong V_{r+1}\op V_{r-1}$ that the
       exterior derivative of a form in the $V_r$ component of $\Lambda^kT^*M$
       has components only in the $V_{r+1}$ and $V_{r-1}$ components of
       $\Lambda^{k+1}T^*M$.  This demonstrates with reassuring simplicity that
       this decomposition naturally leads to a new double complex.  As with the
       Dolbeault complex, the exterior differential $d\omega$ of a $k$-form
       $\omega$ can be split up into the components $\D'\omega$ and $\DD\omega$,
       which are the components of $d\omega$ lying in Sp(1)-representations of
       higher and lower weight respectively.  We describe the operators $\D'$
       and $\DD$ and define new cohomology groups.

In Section \ref{Elipsec} we determine where the (upward) complex is elliptic.
This proves to be a tricky problem, requiring further decomposition and careful
analysis.  Fortunately, this hard work proves to be well worthwhile, as the
double complex is shown to be elliptic in most places.  Like the real Dolbeault
complex, ellipticity only fails at the bottom of the isosceles triangle of
spaces where $\DD=0$ and $\D'=d$, breaking the usual pattern which relies on
studying non-trivial projection maps.

Finally, in Section \ref{Qformsec} we consider some extra opportunities for
decomposition of quaternion-valued forms in the more restricted class of
hypercomplex manifolds.  Here the global complex structures $I$,
$J$ and $K$ can be identified with the action of the imaginary quaternions $i$,
$j$ and $k$ on $\H^n$.  This allows a further decomposition of quaternion-valued
forms, which can be used to develop the theory of quaternionic analysis on
hypercomplex manifolds.  For example, Joyce's algebraic theory of quaternion
holomorphic functions \cite{Jqalg} can be described using this type of
decomposition.


\section{Background Material}
\label{BackSec}

This section contains background information in complex and quaternionic geometry, reviewing the definition of a complex manifold and the ways in which this is adapted to quaternions. We also recall some basic facts concerning the algebra of Sp(1)- 
\linebreak representations.

\subsection{Complex, Hypercomplex and Quaternionic Manifolds}

An {\it almost complex structure} on a $2n$-dimensional real manifold   $M$ is a smooth tensor $I\in C^\infty(\mathrm{End}(TM))$ such that 
$I^2=-\id_{TM}$. 
An almost complex structure $I$ is said to be integrable if and only if the Nijenhuis tensor of $I$
\[N_I(X,Y)=[X,Y]+I[IX,Y]+I[X,IY]-[IX,IY]
\]
vanishes for all $X, Y \in C^\infty(TM)$, for all $x\in M$. 
In this case it can be proved that the almost complex structure $I$ 
arises from a suitable set of holomorphic coordinates on $M$, and 
the pair $(M,I)$ is said to be a complex manifold. 
This way of defining a complex manifold adapts itself well to the quaternions. 

\begin{dfn}
An {\it almost hypercomplex structure} on a $4n$-dimensional manifold $M$ is a triple $(I,J,K)$ of almost complex structures on $M$ which satisfy the relation $IJ=K$. 
An almost hypercomplex structure on $M$ defines an isomorphism 
$T_xM\cong\H^n$ at each point $x\in M$. 

If all of the complex structures are integrable then $(I,J,K)$ is called a {\it hypercomplex structure} on $M$, and $M$ is a 
{\it hypercomplex manifold}.
\label{Hypmfddfn}
\end{dfn}

Not all of the manifolds which we wish to describe as `quaternionic' admit hypercomplex structures. For example, the quaternionic projective line $\HP^1$ is diffeomorphic to the 4-sphere $S^4$. It is well-known that $S^4$ does not even admit a global almost complex structure; so $\HP^1$ can certainly not be hypercomplex, despite behaving extremely like the quaternions locally.

The reason (and the solution) for this difficulty can be 
described succinctly in terms of $G$-structures on manifolds. Let $P$ be the principal frame bundle of $M$, {\it i.e.} the $\GL(n,\R)$-bundle whose fibre over $x\in M$ is the group of isomorphisms $T_xM\cong\R^{4n}$.  
Let $G$ be a Lie subgroup of $\GL(n,\R)$. A {\it $G$-structure} $Q$ on $M$ is a principal subbundle of $P$ with structure group $G$.

Suppose $M^{2n}$ has an almost complex structure. The group of automorphisms of $T_xM$ preserving such a structure is isomorphic to $\GL(n,\C)$. Thus an almost complex structure $I$ and a 
$\GL(n,\C)$-structure $Q$ on $M$ contain the same information. The bundle $Q$ admits a torsion-free connection if and only if there is a
torsion-free linear connection $\nabla$ on $M$ with $\nabla I=0$, in which case it is easy to show that $I$ is integrable. Thus a complex manifold is precisely a real manifold $M^{2n}$ with a GL($n,\C$)-structure $Q$ admitting a torsion-free connection (in which case $Q$ itself is said to be `integrable').

A hypercomplex manifold is a real manifold $M$ with an integrable 
$\GL(n,\H)$-structure $Q$. 
However, the group  
$\GL(n,\H)$ is not the largest subgroup of $\GL(4n,\R)$ preserving the quaternionic structure of $\H^n$. If  
$\GL(n,\H)$ acts on $\H^n$ by right-multiplication by $n\t n$ quaternionic matrices, then the action of $\GL(n,\H)$ commutes with that of the left \H-action of the group 
$\GL(1,\H)\cong \H^*$. 
Thus the group of symmetries of $\H^n$ is the product
$\GL(1,\H)\t_{\R^*}\GL(n,\H)$. Scaling the first factor by a real multiple of the identity reduces the first factor to $\Sp(1)$, and 
$\GL(1,\H)\t_{\R^*}\GL(n,\H)$ is the same as 
$\Sp(1)\t_{\Z_2}\GL(n,\H)$ which is normally abbreviated to 
$\Sp(1)\GL(n,\H)$.
 
\begin{dfn}\cite[1.1]{SalQgeom}
A {\it quaternionic manifold} is a $4n$-dimensional real manifold 
$M$ ($n\ge 2$) with an $\Sp(1)\GL(n,\H)$-structure $Q$ admitting
a torsion-free connection.
\label{Qcmfddfn}
\end{dfn}

When $n=1$ the situation is different, since $\Sp(1)\Sp(1)\cong\mathrm{SO(4)}$. In four dimensions we make the special definition that a quaternionic manifold is a self-dual conformal manifold.

In terms of tensors, quaternionic manifolds are a generalisation of hypercomplex manifolds in the following way. Each tangent space 
$T_xM$ still admits a hypercomplex structure giving an isomorphism 
$T_xM\cong \H^n$, but this isomorphism does not necessarily arise from globally defined complex structures on $M$. There is still an invariant $S^2$-bundle of local almost-complex structures satisfying the equation $IJ=K$, but it is free to `rotate'. 

The manifolds defined above all have Riemannian counterparts defined by a reduction of the structure group $G$ to a compact subgroup. 
A complex manifold $M$ whose $\GL(n,\C)$-structure reduces to an integrable U$(n)$-structure admits a Riemannian metric $g$ with 
$g(IX,IY)=g(X,Y)$ for all $X, Y\in T_xM$ for all $x\in M$. In this case $M$ is called a K\"{a}hler manifold and the metric $g$ is called a K\"{a}hler metric. 
The differentiable 2-form $\omega(X,Y)=g(IX,Y)$, called the 
K\"{a}hler form, is closed and $M$ is a symplectic manifold --- so 
$M$ has compatible complex and symplectic structures.

The quaternionic analogue of the unitary group U$(n)$ is the compact group $\Sp(n)$. A hypercomplex manifold whose $\GL(n,\H)$-structure 
$Q$ reduces to an integrable $\Sp(n)$-structure $Q'$ admits a metric $g$ which is simultaneously K\"{a}hler for each of the complex 
structures $I,J$ and $K$. Such manifolds are called {\it\hk} and have been extensively studied, Fujiki's account \cite{Fujiki} being perhaps the most relevant in this context.
\Hk manifolds have three independent symplectic forms $\omega_I$, $\omega_J$ and 
$\omega_K$. The complex 2-form $\omega_J+i\omega_K$ is holomorphic with respect to the complex structure $I$, and a \hk
manifold has compatible hypercomplex and complex-symplectic structures. 

Similarly, if a quaternionic manifold has a metric compatible with the torsion-free 
\linebreak
$\Sp(1)\GL(n,\H)$-structure, then the 
$\Sp(1)\GL(n,\H)$-structure $Q$ reduces to an $\Sp(1)\Sp(n)$-
\linebreak
structure $Q'$ and $M$ is said to be {\it \qc K\"{a}hler}. 
The group $\Sp(1)\Sp(n)$ is a maximal proper subgroup of SO$(4n)$ except when $n=1$, where $\Sp(1)\Sp(1)=$\ SO(4). In four dimensions a manifold is said to be 
\qc K\"{a}hler if and only if it is self-dual and Einstein.

\subsection{Sp(1)-Representations}

Let Sp(1) be the multiplicative group of unit quaternions. Its Lie algebra $\sp(1)$ is generated by $I$, $J$ and $K$ and the bracket relations $[I,J]=2K$, $[J,K]=2I$ and $[K,I]=2J$. 
The representations of $\Sp(1)$ are ubiquitous in modern mathematics, often in the guise of representations of the isomorphic group $\mathrm{SU}(2)$ or the complexification $\mathrm{SL}(2,\C)$. Standard texts include \cite[\S 2.5]{BD} and \cite[Lecture 11]{Reps}. The second of these is particularly useful for describing Sp(1)-representations and their tensor products using weights of the action of $\sl(2,\C)$. This technique is instrumental in decomposing exterior forms on quaternionic manifolds.  We recall the most salient points.
  
Following standard practice we work primarily with representations on complex vector spaces, real (and quaternionic) representations being constructed in the presence of suitable structure maps. 
Every representation of Sp(1) on a complex vector space can be written as a 
direct sum of irreducible representations and the multiplicity of each irreducible in such a decomposition is uniquely determined. 

Let $V_1$ be the basic representation of $\Sp(1)\cong\mathrm{SU}(2)$ on $\C^2$ given by left-action of matrices upon column vectors. 
The $n^{\mathrm{th}}$ symmetric power of $V_1$ is a   
representation on $\C^{n+1}$ which is written   
\[ V_n=S^n(V_1).
\] 
The representation $V_n$ is irreducible  
and every irreducible representation of Sp(1) is of the form $V_n$ for some nonnegative $n\in\Z$.  
Each irreducible representation $V_n$ is an eigenspace of the Casimir operator $I^2+J^2+K^2$ with eigenvalue $-n(n+2)$. 

The irreducible representation $V_n$ can be decomposed into weight spaces under the action of a Cartan subalgebra of $\sl(2,\C)$. Each weight space is one-dimensional and the weights are the integers 
\[\{n,n-2,\ldots,n-2k,\ldots,2-n,-n\}.
\]
Thus $V_n$ is also characterised by being the unique irreducible representation of $\sl(2,\C)$ with highest weight $n$. 

It follows from the Leibniz rule $I(a\ot b)=I(a)\ot b + a\ot I(b)$ for Lie algebra representations that if $a\in V_m$ and $b\in V_n$ are weight vectors of $I$ with weights $\lambda$ and $\mu$ then $a\ot b$ is a weight vector of 
$V_m\ot V_n$ with weight $\lambda+\mu$.  
Weight space decompositions can thus be used to determine tensor, symmetric and exterior products of Sp(1)-representations. 
Amongst other things, this enables us to calculate the irreducible decomposition of the (diagonal) action of 
Sp(1) on the tensor product 
$V_m\ot V_n$. This is given by the famous {\it Clebsch-Gordon formula},
\be \quad\quad\ \ \ \ \ \ \ \ \ \ \ \ V_m\ot V_n
\cong V_{m+n}\op V_{m+n-2}\op\cdots\op V_{m-n+2}\op V_{m-n}\quad
\ \ \ \mbox{for $m\geq n$}.
\label{CGform}
\ee

\subsection{Associated Bundles}

It is usual in differential geometry to talk about properties of a vector bundle using the properties of its fibres. A good example of this is when talking about a group acting on a vector bundle: or more precisely, a principal bundle whose fibres act upon the fibres of an associated vector bundle. 

Let $P$ be a principal 
$G$-bundle over the differentiable manifold $M$ and let $V$ be a representation of the group $G$. We define the {\it associated bundle}  
\[\mathbf{V}=P\t_G V= \frac{P\t V}{G},\]
where $g\in G$ acts on $(p,v)\in P\t V$ by 
$(p,v)\cdot g=(f\cdot g,g^{-1}\cdot v)$. Then $\mathbf{V}$ is a vector bundle over $M$ with fibre $V$. At every point $m\in M$ the fibre $P_m\cong G$ acts on the fibre $\mathbf{V}_m$, a notion that is commonly abused slightly by saying that $G$ acts on $\mathbf{V}$. 
A decomposition of $V$ into subrepresentations then gives rise to a decomposition of 
$\mathbf{V}$ into associated subbundles. 
We will usually just write 
$V$ for $\mathbf{V}$, relying on context to distinguish between the bundle and the representation. 

In this way the cotangent spaces of complex and quaternionic manifolds are described using representations of the groups 
$\GL(n,\C)$ and $\Sp(1)\GL(n,\H)$ respectively: or indeed using the much simpler representations of the groups U(1) and $\Sp(1)$.


\section{Differential Forms on Complex Manifolds}
\label{Dolsec}

Let $(M,I)$ be a complex manifold. 
The complexified cotangent space splits into eigenspaces of $I$ with eigenvalues $\pm i$, $T^*M^\C=T^*_{1,0}M\oplus T^*_{0,1}M$, which are called the holomorphic and antiholomorphic cotangent spaces respectively.  
This induces the familiar decomposition into types of exterior 
$k$-forms  
\[ \Lambda^kT^*M^\C=\bigoplus_{p+q=k}\Lambda^p(T^*_{1,0}M)\ot
                                    \Lambda^q(T^*_{0,1}M), 
\]
where the bundle
$\Lambda^{p,q}M=\Lambda^pT^*_{1,0}M\ot\Lambda^qT^*_{0,1}M$ 
is called the bundle of $(p,q)$-forms on $M$.
A smooth section of the bundle $\Lambda^{p,q}M$ is called a  differential form of type $(p,q)$ or just a $(p,q)$-form. We write $\Omega^{p,q}(M)$ for the set of $(p,q)$-forms on $M$, so
\[\Omega^{p,q}(M)=C^\infty(M,\Lambda^{p,q}M)\quad\quad\mbox{and}
   \quad\quad\Omega^k(M)=\bigoplus_{p+q=k}\Omega^{p,q}(M).\]
Define two first-order differential operators,
\be 
\begin{array}{l}
\partial:\Omega^{p,q}(M)\ra \Omega^{p+1,q}(M)\\ 
    \partial=\pi^{p+1,q}\circ d
\end{array}
\quad\mathrm{and}\quad
\begin{array}{l}
\Dol: \Omega^{p,q}(M)\ra \Omega^{p,q+1}(M)\\
    \Dol=\pi^{p,q+1}\circ d,
\end{array}
\label{Doldfn}
\ee
where $\pi^{p,q}$ denotes the natural projection from 
$\Lambda^kT^*M^\C$ onto 
$\Lambda^{p,q}M$.
The operators $\partial$ and $\Dol$ are called the Dolbeault operators.

These definitions rely only on the fact that $I$ is an almost complex structure. If in addition $I$ is integrable then these are the only two components of  
the exterior differential $d$, so that $d=\partial+\Dol$ 
\cite[p.34]{W}. An immediate consequence of this is that on a complex manifold $M$, 
$\partial^2=\partial\Dol+\Dol\partial=\Dol^2=0$.
This gives rise to the {\it Dolbeault complex}. 

\begin{figure}
\caption{The Dolbeault Complex}
\begin{center}
\begin{picture}(300,210)

\put(0,100){\makebox(20,0)[c]{$\Omega^{0,0}=C^\infty(M,\C)$}}
\put(20,110){\vector(2,1){60}}
\put(20,90){\vector(2,-1){60}}
\put(50,130){\makebox(0,0)[b]{$\partial$}}
\put(50,70){\makebox(0,0)[t]{$\Dol$}}

\put(100,150){\makebox(0,0)[c]{$\Omega^{1,0}$}}
\put(120,160){\vector(2,1){60}}
\put(120,140){\vector(2,-1){60}}
\put(150,180){\makebox(0,0)[b]{$\partial$}}
\put(150,120){\makebox(0,0)[t]{$\Dol$}}

\put(100,50){\makebox(0,0)[c]{$\Omega^{0,1}$}}
\put(120,60){\vector(2,1){60}}
\put(120,40){\vector(2,-1){60}}
\put(150,80){\makebox(0,0)[b]{$\partial$}}
\put(150,20){\makebox(0,0)[t]{$\Dol$}}

\put(200,100){\makebox(0,0)[c]{$\Omega^{1,1}$}}
\put(220,110){\vector(2,1){60}}
\put(220,90){\vector(2,-1){60}}
\put(250,130){\makebox(0,0)[b]{$\partial$}}
\put(250,70){\makebox(0,0)[t]{$\Dol$}}

\put(200,200){\makebox(0,0)[c]{$\Omega^{2,0}$}}
\put(220,190){\vector(2,-1){60}}
\put(250,170){\makebox(0,0)[t]{$\Dol$}}
\put(250,200){\makebox(0,0)[l]{$\ldots\ldots$ {\it etc.}}}

\put(200,0){\makebox(0,0)[c]{$\Omega^{0,2}$}}
\put(220,10){\vector(2,1){60}}
\put(250,30){\makebox(0,0)[b]{$\partial$}}
\put(250,0){\makebox(0,0)[l]{$\ldots\ldots$ {\it etc.}}}

\end{picture}
\end{center}
\end{figure}
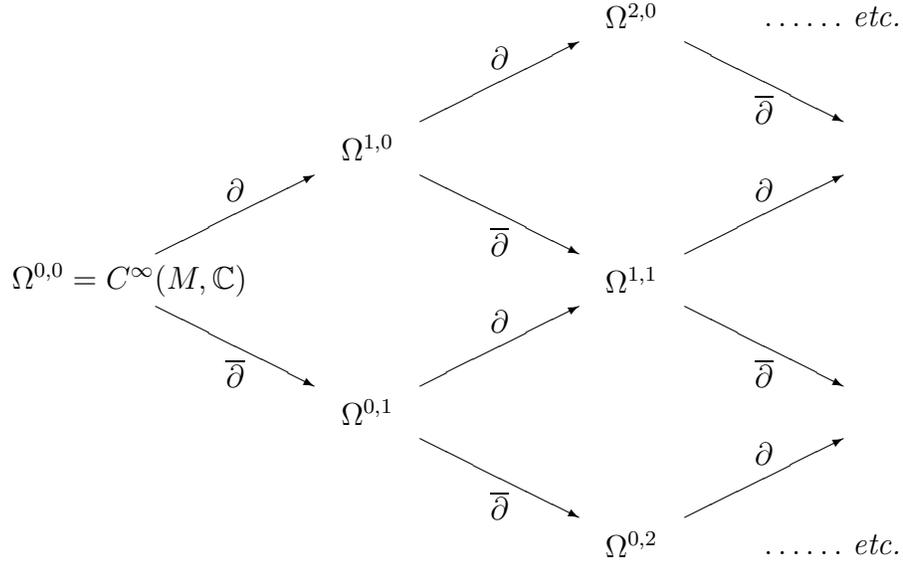

The purpose of this paper is precisely to present the quaternionic analogue of this double complex. To do this, 
note that the bundle $\Lambda^{p,q}M$ is an eigenspace of the induced action of $I$ on 
$\Lambda^kT^*M^\C$, since for $\omega\in \Lambda^{p,q}M$, 
$I(\omega)=i(p-q)\omega$.   
The decomposition into types can therefore be thought of as a decomposition of $\Lambda^kT^*M^\C$ into $\u(1)$-representations, where the complex structure $I$ generates a copy of the Lie algebra $\u(1)$. 
\footnote{The irreducible representations of U(1) are all one-dimensional. They are parametrised by the integers, taking the form 
$\varrho_n:\mathrm{U}(1)\ra \GL(1,\C)=\C^*$, 
$\varrho_n:e^{i\theta}\ra e^{ni\theta}$ 
for some $n\in\Z$. The corresponding representations of the Lie algebra 
${\goth u}(1)$ are then of the form $d\varrho_n:z\ra inz$.} 
We will presently do this in quaternionic geometry by replacing the Lie algebra 
$\u(1)=\langle I \rangle$ with $\sp(1)=\langle I,J,K \rangle$ and decompose $\Lambda^kT^*M$ into irreducible $\sp(1)$-representations.

\subsection{Real forms on Complex Manifolds}

It is less well-known that a similar splitting occurs for real-valued exterior forms. This is an instructive case, because the resulting double complex is even more closely akin to the new quaternionic double complex.
 
Let $M$ be a complex manifold and let 
$\omega\in \Lambda^{p,q}=\Lambda^{p,q}M$. For simplicity's sake assume that $p>q$ throughout.  
Then $\ov{\omega}\in\Lambda^{q,p}$, and $\omega+\ov{\omega}$ is a real-valued exterior form. 
Define the space of such forms, 
\[	(\Lambda^{p,q}\op\Lambda^{q,p})_\R
				=[\Lambda^{p,q}\op\Lambda^{q,p}]
				\equiv[[\Lambda^{p,q}]].
\]
The space    
$[[\Lambda^{p,q}]]$ is a real vector bundle associated to the principal  
GL$(n,\C)$-bundle defined by the complex structure. The first square bracket indicates real forms and the second the direct sum, the  notation following that of Reyes-Carri\'{o}n  \cite[\S3.1]{R-C}, who uses the ensuing decomposition on K\"{a}hler manifolds.   

This gives a decomposition of real-valued exterior forms,
\be \Lambda_\R^kT^*M=\bigoplus_{\stackrel{\scriptstyle{p+q=k}}{p>q}}
                 [[\Lambda^{p,q}]]
                 \op[\Lambda^{\frac{k}{2},\frac{k}{2}}].
\label{realDol}
\ee
The condition $p>q$ ensures that we have no repetition. 
The bundle $[\Lambda^{\frac{k}{2},\frac{k}{2}}]$ only appears when $k$ is even. It is its own conjugate and so naturally a real vector bundle associated to the trivial representation ({\it i.e.} the zero weight space) of the Lie algebra 
$\u(1)=\langle I \rangle$. 

Let $[[\Omega^{p,q}]]=C^\infty([[\Lambda^{p,q}]])$, so that for   
$\omega\in\Omega^{p,q}(M)$, 
$\omega+\ov{\omega}\in [[\Omega^{p,q}]]$. Then 
\[
d(\omega+\ov{\omega}) = (\partial\omega+\Dol\ov{\omega})+
                        (\Dol\omega+\partial\ov{\omega}) 
\in[[\Omega^{p+1,q}]]\op[[\Omega^{p,q+1}]].
\]
Call the first of these components $[\partial]\omega$ and the second $[\Dol]\omega$. This defines real analogues of the Dolbeault operators. The `double complex' equations 
$[\partial]^2=[\partial][\Dol]+[\Dol][\partial]=[\Dol]^2=0$ 
follow directly from decomposing the equation $d^2=0$.  
For the space $[\Lambda^{\frac{k}{2},\frac{k}{2}}]$ there is no space $[[\Lambda^{\frac{k}{2},\frac{k}{2}+1}]]$ because $p>q$, so 
$[\Dol]=0$ and there is just one operator $[\partial]=d$.

\begin{figure}[h]
\caption{The Real Dolbeault Complex}
\begin{center}
\begin{picture}(400,160)

\put(0,0){\makebox(20,0)[c]{$[\Omega^{0,0}]=C^\infty(M)$}}
\put(20,10){\vector(2,1){60}}
\put(45,27){\makebox(0,0)[b]{$d$}}

\put(100,50){\makebox(0,0)[c]{$[[\Omega^{1,0}]]=\Omega^1(M)$}}
\put(120,60){\vector(2,1){60}}
\put(120,40){\vector(2,-1){60}}
\put(135,77){\makebox(0,0)[b]{$[\partial]$}}
\put(135,23){\makebox(0,0)[t]{$[\Dol]$}}

\put(200,0){\makebox(0,0)[c]{$[\Omega^{1,1}]$}}
\put(220,10){\vector(2,1){60}}
\put(245,27){\makebox(0,0)[b]{$d$}}

\put(200,100){\makebox(0,0)[c]{$[[\Omega^{2,0}]]$}}
\put(220,110){\vector(2,1){60}}
\put(220,90){\vector(2,-1){60}}
\put(235,127){\makebox(0,0)[b]{$[\partial]$}}
\put(235,73){\makebox(0,0)[t]{$[\Dol]$}}

\put(300,50){\makebox(0,0)[c]{$[[\Omega^{2,1}]]$}}
\put(320,60){\vector(2,1){60}}
\put(320,40){\vector(2,-1){60}}
\put(335,77){\makebox(0,0)[b]{$[\partial]$}}
\put(335,23){\makebox(0,0)[t]{$[\Dol]$}}
\put(350,50){\ldots\ldots{\it etc.}}

\put(300,150){\makebox(0,0)[c]{$[[\Omega^{3,0}]]$}}
\put(320,140){\vector(2,-1){60}}
\put(335,123){\makebox(0,0)[t]{$[\Dol]$}}
\put(350,150){\ldots\ldots{\it etc.}}

\end{picture}
\end{center}
\label{realdol}
\end{figure}
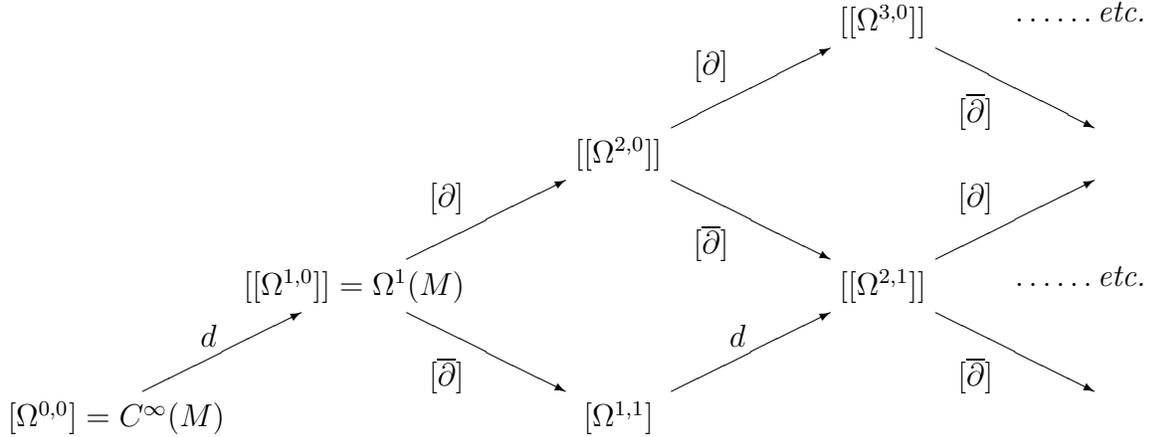

Thus there is a double complex of real forms on a complex manifold, obtained by decomposing $\Lambda^k_xT^*M$ into subrepresentations of the action of ${\goth u}(1)=\langle I\rangle$, induced from the action on $T^*_xM$. 
The main difference is that this real-valued complex gives an isosceles triangle of spaces, whereas the standard Dolbeault complex gives a full diamond. Amalgamating the spaces $\Lambda^{p,q}$ and 
$\Lambda^{q,p}$ into the single real space $[[\Lambda^{p,q}]]$ has effectively folded this diamond in half. This structure is very similar to that of the new quaternionic double complex which is the main subject of this paper.


\subsubsection*{Ellipticity}

A complex $0\stackrel{\Phi_0}{\lra}C^\infty(E_0)
            \stackrel{\Phi_1}{\lra}C^\infty(E_1)
            \stackrel{\Phi_2}{\lra}C^\infty(E_2)
            \stackrel{\Phi_3}{\lra}\ldots
            \stackrel{\Phi_n}{\lra}C^\infty(E_n) 
            \stackrel{\Phi_{n+1}}{\lra}0$
is said to be elliptic at $E_i$ if the 
{\it principal symbol sequence}  
$E_{i-1}\stackrel{\sigma_{\Phi_i}}{\lra}E_i
  \stackrel{\sigma_{\Phi_{i+1}}}{\lra}E_{i+1}$
 is exact for all 
$\xi\in T_x^*M$ and for all $x\in M$. 
\footnote{The link between elliptic complexes and elliptic operators (those whose principal symbol is an isomorphism) is as follows. Given any metric on each
$E_i$, define a formal adjoint $\Phi^*_i:E_i\ra E_{i-1}$.
The complex is elliptic at $E_i$ if and only if the Laplacian $\Phi^*_i\Phi_i+\Phi_{i-1}\Phi_{i-1}^*$ is an elliptic operator.} 
A thorough description of this topic can be found in 
\cite[Chapter 5]{W}. For our operators, it suffices to note that 
principal symbol of the exterior differential 
$d:\Omega^k(M)\ra\Omega^{k+1}(M)$ is 
$\sigma_d(x,\xi)\omega=\omega\wedge\xi$. Let $\pi:\Lambda^kT^*M\ra E$ be the projection from $\Lambda^k T^*M$ onto some subspace 
$E$. It follows that the principal symbol of $\pi\circ d$ on $\Lambda^{k-1}T^*M$ is just  
$\sigma_{\pi\circ d}(x,\xi)\omega=\pi(\omega\wedge\xi)$.  

It is important to establish where a differential complex is elliptic for various reasons: for example, an elliptic complex on a compact manifold always has finite-dimensional cohomology groups 
\cite[Theorem 5.2, p. 147]{W}. 
The de Rham 
and Dolbeault complexes are elliptic everywhere. 
The real Dolbeault complex of Figure \ref{realdol} is elliptic in most places, but not everywhere.  Interestingly, it fails to be elliptic in almost exactly the same places as the new quaternionic double complex, and for the same reasons.

\begin{prop}
For $p>0$, the upward complex 
\[0\lra [\Omega^{p,p}]\stackrel{d}{\lra} [[\Omega^{p+1,p}]]
\stackrel{[\partial]}{\lra} [[\Omega^{p+2,p}]]\stackrel{[\partial]}{\lra}\ldots
\]
is elliptic everywhere except at the first two spaces 
$[\Omega^{p,p}]$ and 
$[[\Omega^{p+1,p}]]$. 

For $p=0$, the `leading edge' complex 
\[0\lra [\Omega^{0,0}]\stackrel{d}{\lra}  [[\Omega^{1,0}]]
\stackrel{[\partial]}{\lra} [[\Omega^{2,0}]]\stackrel{[\partial]}{\lra}\ldots
\]
is elliptic everywhere except at 
$[[\Omega^{1,0}]]=\Omega^1(M)$.
\label{RealDolelipthm}
\end{prop}

\begin{proof}
When $p>q+1$, the short sequence 
$[[\Omega^{p-1,q}]]
\stackrel{[\partial]}{\lra} [[\Omega^{p,q}]]\stackrel{[\partial]}{\lra}
[[\Omega^{p+1,q}]]$
is a real form of the sequence 
\[
\begin{array}{ccccc}
\Omega^{p-1,q} & \stackrel{\partial}{\lra} & \Omega^{p,q} &
                 \stackrel{\partial}{\lra} & \Omega^{p+1,q} \\
\bigoplus &  & \bigoplus &  & \bigoplus \\
\Omega^{q,p-1} & \stackrel{\Dol}{\lra} & \Omega^{q,p} &
                 \stackrel{\Dol}{\lra} & \Omega^{q,p+1}.
\end{array}
\]
This is (a real subspace of) the direct sum of two elliptic sequences, and so is elliptic. Thus we have ellipticity at 
$[[\Omega^{p,q}]]$ whenever $p\geq q+2$. 

This leaves us to consider the case when $p=q$, giving (a real subspace of) the sequence 
\be
\begin{array}{ccccc}
0 & \lra &\Omega^{p,p}& 
\begin{array}{c}\stackrel{\partial}{\nearrow} \\ 
                \stackrel{\textstyle{\searrow}}{\scriptstyle{\Dol}}
\end{array}
&
\begin{array}{cccc}\Omega^{p+1,p} &
    \stackrel{\partial}{\lra} & \Omega^{p+2,2} 
                    &\lra\ldots\mbox{{\it etc.}} \\
    \bigoplus & &\bigoplus                    \\
    \Omega^{p,p+1} &
                 \stackrel{\Dol}{\lra} & \Omega^{2,p+2}
                    &\lra\ldots\mbox{{\it etc.}} \\
\end{array}
\end{array}
\ee
This fails to be elliptic. 
An easy and instructive way to see this is to consider the simplest 4-dimensional example $M=\C^2$.

Let $e^0$, $e^1=I(e^0)$, $e^2$ and $e^3=I(e^2)$ form a basis for 
$T^*_x\C^2\cong\C^2$, and let $e^{ab\ldots}$ denote 
$e^a\wedge e^b\wedge\ldots\mbox{\it etc.}$ Then 
$I(e^{01})=e^{00}-e^{11}=0$, so $e^{01}\in[\Lambda^{1,1}]$. 
The map from $[\Lambda^{1,1}]$ to $[[\Lambda^{2,1}]]$ is just the exterior differential $d$. Since  
$\sigma_d(x,e^0)(e^{01})=e^{01}\wedge e^0=0$ the symbol map $\sigma_d:[\Lambda^{1,1}]\ra[[\Lambda^{2,1}]]$ is not injective, so the symbol sequence is not exact at $[\Lambda^{1,1}]$. 

Consider also $e^{123}\in[[\Lambda^{2,1}]]$. Then  
$\sigma_{[\partial]}(x,e^0)(e^{123})=0$, since there is no bundle 
$[[\Lambda^{3,1}]]$. But $e^{123}$ has no $e^0$-factor, so is not the image under $\sigma_d(x,e^0)$ of any form 
$\alpha\in[\Lambda^{1,1}]$. Thus the symbol sequence fails to be exact at $[[\Lambda^{2,1}]]$.   
 
It is a simple matter to extend these counterexamples to higher dimensions and higher exterior powers. For $k=0$, the  
situation is different. It is easy to show that the complex 
\[0\lra C^\infty(M)\stackrel{d}{\lra}
       [[\Omega^{1,0}]]
  \stackrel{[\partial]}{\lra}[[\Omega^{2,0}]] 
  \lra\ldots\mbox{{\it etc.}}
\]
is elliptic everywhere except at $[[\Omega^{1,0}]]$. 
\end{proof}

This last sequence is given particular attention by 
Reyes-Carri\'{o}n \cite[Lemma 2]{R-C}. He shows that, when $M$ is 
K\"{a}hler, ellipticity can be regained by adding the space 
$\langle\omega\rangle$ to the bundle 
$[[\Lambda^{2,0}]]$, 
where $\omega$ is the real K\"{a}hler $(1,1)$-form.

The real Dolbeault complex is thus elliptic except at the bottom of the isosceles triangle of spaces. Here the projection from 
$d([\Omega^{p,p}])$ to 
$[[\Omega^{p+1,p}]]$ is the identity, and arguments based upon non-trivial projection maps no longer apply. 
We shall see that this situation is closely akin to that of differential forms on quaternionic manifolds, and that techniques motivated by this example yield similar results.


\section{Differential Forms on Quaternionic Manifolds}
\label{Backgroundsec}

This section describes previous results in the decomposition of exterior forms in quaternionic geometry. These fall into two categories: those arising from taking repeated products with the fundamental 4-form in quaternionic K\"{a}hler geometry and those arising from considering the representations of $\GL(n,\H)\Sp(1)$ on $\Lambda^kT^*M$. We are primarily concerned with the second approach.

The decomposition of differential forms on quaternionic K\"{a}hler manifolds began by considering the 
{\it fundamental $4$-form}
\[\Omega=\omega_I\wedge\omega_I+\omega_J\wedge\omega_J
        +\omega_K\wedge\omega_K,
\]
where $\omega_I$, $\omega_J$ and $\omega_K$ are the local K\"{a}hler forms associated to local almost complex structures $I, J$ and $K$  with $\,IJ=K$. 
The fundamental 4-form is globally defined and invariant under the induced action of Sp(1)Sp(n) on $\Lambda^4T^*M$. 
Kraines \cite{YK} and Bonan \cite{Bonan} used the fundamental 4-form to decompose the space $\Lambda^kT^*M$ in a similar way to the 
Lefschetz decomposition of differential forms on a K\"{a}hler manifold 
\cite[p. 122]{GH}. 
A differential $k$-form $\mu$ is said to be 
{\it effective} if $\;\Omega\wedge*\mu=0$, where 
$*:\Lambda^kT^*M\ra\Lambda^{4n-k}T^*M$ is the Hodge star. This leads to the following theorem: 

\begin{thm}
\cite[Theorem 3.5]{YK}\cite[Theorem 2]{Bonan}
Let $M^{4n}$ be a quaternionic K\"{a}hler manifold. 
For $k\leq 2n+2$, every every $k$-form $\phi$ admits a unique decomposition
\[\phi=\sum_{0\leq j\leq k/4}\Omega^j\wedge\mu_{k-4j},\]
where the $\mu_{k-4j}$ are effective $(k-4j)$-forms.
\label{Bonandecomp}
\end{thm}

Bonan further refines this decomposition for quaternion-valued forms, using exterior multiplication by the globally defined quaternionic 2-form $\Psi=i_1\omega_I+i_2\omega_J+i_3\omega_K$. Note that $\Psi\wedge\Psi=-2\Omega$.

Another way to consider the decomposition of forms on a quaternionic manifold is as subbundles of $\Lambda^kT^*M$ associated with different representations of the group \linebreak
$\Sp(1)\GL(n,\H)$. 
The representation of $\Sp(1)\GL(n,\H)$ 
 on $\H^n$ is given by the equation 
\be\H^n\ot_\R\C\cong V_1\ot E,
\label{Hnrepeqn}
\ee
where $V_1$ is the basic representation of Sp(1) on $\C^2$ and 
$E$ is the basic representation of $\,\GL(n,\H)$\, on $\C^{2n}$. (This uses the standard convention of working with complex  representations, which in the presence of suitable structure maps can be thought of as complexified real representations. In this case, the structure map is the tensor product of the quaternionic structures on $V_1$ and $E$.)

This representation also describes the (co)tangent bundle of a quaternionic manifold in the following way. 
Following Salamon \cite[\S1]{SalQgeom}, if $M^{4n}$ is a quaternionic manifold with $\Sp(1)\GL(n,\H)$-structure $Q$, then the cotangent   
bundle is a vector bundle associated with the principal bundle $Q$ and the representation $V_1\ot E$, so that  
\be (T^*M)^\C\cong V_1\ot E
\label{TMrepeqn}
\ee
(though we will usually omit the complexification sign). 
This induces an $\Sp(1)\GL(n,\H)$-action on the bundle of exterior 
$k$-forms $\Lambda^k T^*M$,
\be\Lambda^kT^*M\cong\Lambda^k(V_1\ot E)\cong 
                 \bigoplus_{j=0}^{[k/2]}S^{k-2j}(V_1)\ot L^k_j
\cong \bigoplus_{j=0}^{[k/2]}V_{k-2j}\ot L^k_j,
\label{Swanndecomp}
\ee
where $L^k_j$ is an irreducible representation of \,GL$(n,\H)$.
This decomposition is given by Salamon \cite[\S 4]{SalQgeom}, along with more details concerning the nature of the $\GL(n,\H)$ representations $L^k_j$. 

If the $\Sp(1)\GL(n,\H)$-structure on $M$ reduces to an 
$\Sp(1)\Sp(n)$-structure, 
$\Lambda^kT^*M$ can be further decomposed into representations of the compact group $\Sp(1)\Sp(n)$. This refinement is performed in detail by Swann \cite{Swann}, 
and used to demonstrate that if $\dim M\ge 8$, the vanishing condition $\nabla\Omega=0$ implies that $M$ is quaternionic 
K\"{a}hler for any torsion-free connection $\nabla$ preserving the $\Sp(1)\Sp(n)$-structure. 
 
If we symmetrise completely on $V_1$ in Equation (\ref{Swanndecomp}) to obtain $V_k$, we must antisymmetrise completely on $E$. Salamon therefore defines the irreducible subspace 
\be A^k\cong V_k\ot\Lambda^k E.
\label{Akdfn}
\ee
The bundle $A^k$ can be described using the decomposition into types for the local almost complex structures on $M$ as follows \cite[Proposition 4.2]{SalQgeom}:
\footnote{This is because every Sp(1)-representation $V_n$ is generated by its highest weight spaces taken with respect to all the different linear combinations of $I$, $J$ and $K$.}
\be A^k=\sum_{I\in S^2}\Lambda_I^{k,0}M.
\label{Salgenthm}
\ee
Letting $p$ denote the natural projection $p:\Lambda^kT^*M\ra A^k$ and setting $D=p\circ d$, 
Salamon defines a sequence of differential operators 
\be 0\longrightarrow C^\infty(A^0)\stackrel{D=d}{\longrightarrow}
   C^\infty(A^1=T^*M)\stackrel{D}{\longrightarrow}
   C^\infty(A^2)\stackrel{D}{\longrightarrow}
   \ldots
   \stackrel{D}{\longrightarrow} C^\infty(A^{2n})\longrightarrow 0.
\label{Salcomplex}
\ee
This is accomplished using only the fact that $M$ has an 
$\Sp(1)\GL(n,\H)$-structure; such a manifold is called `almost quaternionic'.
The following theorem of Salamon relates the integrability of such a structure with the sequence of operators in (\ref{Salcomplex}):

\begin{thm}\cite[Theorem 4.1]{SalQgeom}
An almost quaternionic manifold is quaternionic if and only if (\ref{Salcomplex}) is a complex.
\label{Salcomplexthm}
\end{thm}

This theorem is analogous to the familiar result in complex geometry that an almost complex structure on a manifold is integrable if and only if $\Dol^2=0.$


\section{Construction of the Double Complex}
\label{DCsec}

In this, the most important section of this paper, we construct the new double complex on a quaternionic manifold $M$ by decomposing the action of $\Sp(1)$ on $\Lambda^kT^*M$ inherited from the 
$\Sp(1)\GL(n,\H)$-structure. The top row of this double complex is the complex (\ref{Salcomplex}) discovered by Salamon.

Let $M^{4n}$ be a quaternionic manifold. Following Salamon 
\cite[\S1]{SalQgeom} we can define (at least locally) vector bundles $\mathbf{V_1}$ and $\mathbf{E}$ associated to the basic complex representations of Sp(1) and 
$\GL(n,\H)$ respectively, so that 
$T^*_xM\cong (\mathbf{V_1})_x\ot\mathbf{E_x}\cong V_1\ot E$ as an $\Sp(1)\GL(n,\H)$-representation for all $x \in M$. Suppose we consider just the action of the Sp(1)-factor. Then the (complexified) cotangent space effectively takes the form 
$V_1\ot\C^{2n}\cong 2nV_1$. 
Whilst the bundles $\mathbf{V_1}$ and $\mathbf{E}$ might be neither globally nor uniquely defined, the Casimir operator $I^2+J^2+K^2$ 
is invariant. It follows that, though the Sp(1)-action on 
a $k$-form $\alpha$ might be subject to choice, its spectrum under the Casimir action, and hence its decomposition into Sp(1)-representations of different weights, is uniquely and globally defined by the $\Sp(1)\GL(n,\H)$-structure. Thus the irreducible decomposition of the Sp(1)-action on 
$\Lambda^kT^*_xM$ is given by the irreducible decomposition of the representation $\Lambda^k(2nV_1)$. 

To work out the irreducible decomposition of this representation 
we compute the weight space decomposition of 
$\Lambda^k(2nV_1)$ from that of $2nV_1$. 
\footnote{This process for calculating the weights of tensor, symmetric and exterior powers is a standard technique in representation theory --- see for example \cite[\S11.2]{Reps}.}
With respect to the action of a particular subgroup 
U$(1)\subset\Sp(1)$,
the representation $2nV_1$ has weights $+1$ and $-1$, each occuring with multiplicity $2n$. 
The weights of $\Lambda^k(2nV_1)$ are the $k$-wise distinct sums of these.  
Each weight $r$ in $\Lambda^k(2nV_1)$ must therefore be a sum of $p$ occurences of the weight `+1' and $p-r$ occurences of the weight 
`$-1$', where $2p-r=k$ and
$0\leq p\leq k$ (from which it follows immediately that 
$-k\leq r\leq k$ and $r\equiv k\ \mathrm{mod}\ 2$).  
The number of ways to choose the $p$ `+1' weights is
$2n\choose p$, 
and the number of ways to choose the $(p-r)$ `$-1$' weights is
$2n\choose p-r$, 
so the multiplicity of the weight $r$ in the representation 
$\Lambda^k(2nV_1)$ is 
\[\mathrm{Mult}(r)={2n\choose{\frac{k+r}{2}}} {2n\choose{\frac{k-r}{2}}}.
\]

For $r\geq 0$, consider the difference Mult$(r)-$Mult$(r+2)$. This is the number of weight spaces of weight $r$ which do not have any corresponding weight space of weight $r+2$. Each such weight space must therefore be the highest weight space in an irreducible subrepresentation $V_r\subseteq\Lambda^kT^*M$, from which it follows that the number of irreducibles $V_r$ in $\Lambda^k(2nV_1)$ is equal to 
Mult$(r)-$Mult$(r+2)$. This demonstrates the following proposition:

\begin{prop}Let $M^{4n}$ be a hypercomplex manifold. The decomposition into irreducibles of the induced representation of 
$\Sp(1)$ on $\Lambda^kT^*M$ is 
\[\Lambda^kT^*M\cong\bigoplus_{r=0}^k
       \left[{2n\choose\frac{k+r}{2}}{2n\choose\frac{k-r}{2}}-
       {2n\choose\frac{k+r+2}{2}}{2n\choose\frac{k-r-2}{2}}\right]V_r,
\]
where $r\equiv k \bmod 2$.
\label{LambdaMdecomp}
\end{prop}

We will not always write the condition $r\equiv k \bmod 2$, assuming that ${p\choose q}=0$ if $q\not\in\Z$.

\begin{dfn}
Let $M^{4n}$ be a quaternionic manifold.
Define $E_{k,r}$ to be the vector subbundle of $\Lambda^kT^*M$ consisting of Sp(1)-representations with highest weight $r$. 
Define the coefficient 
\[ \epsilon^n_{k,r}= 
{2n\choose\frac{k+r}{2}}{2n\choose\frac{k-r}{2}}-
       {2n\choose\frac{k+r+2}{2}}{2n\choose\frac{k-r-2}{2}},
\]
so that (neglecting the $\GL(n,\H)$-action) we have 
$E_{k,r}\cong\epsilon^n_{k,r}V_r$. 
\label{Enotationdfn}
\end{dfn}

With this notation Proposition \ref{LambdaMdecomp} is written
\[\Lambda^kT^*M\cong \bigoplus_{r=0}^k \epsilon^n_{k,r}V_r
               \cong \bigoplus_{r=0}^k E_{k,r}.
\]
Our most important result is that this decomposition gives rise to a
double complex of differential forms and operators on a quaternionic
manifold.

\begin{thm}
The exterior derivative $d$ maps $C^\infty(M,E_{k,r})$ to 
$C^\infty(M,E_{k+1,r+1}\oplus E_{k+1,r-1})$.
\label{1stDCthm}
\end{thm}

\begin{proof}
Let $\nabla$ be a torsion-free linear
 connection on $M$ preserving the quaternionic structure.
Then $\nabla:C^\infty(M,E_{k,r})\ra 
                     C^\infty(M,E_{k,r}\ot T^*M)$. 
As Sp(1)-representations, this is
\[\nabla:C^\infty(M,\epsilon^n_{k,r}V_r)\ra 
C^\infty(M,\epsilon^n_{k,r}V_r\ot 
2nV_1).
\]
Using the Clebsch-Gordon formula we have 
$\epsilon^n_{k,r}V_r\ot 2nV_1\cong 
2n\epsilon^n_{k,r}(V_{r+1}\op V_{r-1})$. 
Thus
the image of $E_{k,r}$ under $\nabla$ is contained in the 
$V_{r+1}$ and $V_{r-1}$ summands of $\Lambda^kT^*M\ot T^*M$.
Since $\nabla$ is torsion-free, 
$d=\wedge\circ\nabla$, so $d$ maps (sections of) $E_{k,r}$ 
to the $V_{r+1}$ and $V_{r-1}$ summands of $\Lambda^{k+1}T^*M$.
Thus $d:C^\infty(M,E_{k,r})\ra 
C^\infty(M,E_{k+1,r+1}\oplus E_{k+1,r-1})$.
\end{proof}

This allows us to split the exterior differential $d$ into two 
`quaternionic Dolbeault operators'.

\begin{dfn}
Let $\pi_{k,r}$ be the natural projection from $\Lambda^kT^*M$ onto
$E_{k,r}$.
Define the operators 
\be 
\begin{array}{l}
              \D':C^\infty(E_{k,r})\ra C^\infty(E_{k+1,r+1})\\ 
              \D'=\pi_{k+1,r+1}\circ d
\end{array}
\quad\mbox{and}\quad
\begin{array}{l}
              \DD:C^\infty(E_{k,r})\ra C^\infty(E_{k+1,r-1})\\
              \DD=\pi_{k+1,r-1}\circ d.
\end{array}
\label{QDoldfn}
\ee
\end{dfn}

Theorem \ref{1stDCthm} is equivalent to the following:

\begin{prop}
On a quaternionic manifold $M$, we have $d=\D'+\DD$, and so
\[\D'^2=\D'\DD+\DD\D'=\DD^2=0.\]
\end{prop}

\begin{proof}
The first equation is equivalent to Theorem \ref{1stDCthm}.
The rest follows immediately from decomposing the equation $d^2=0$.
\end{proof}


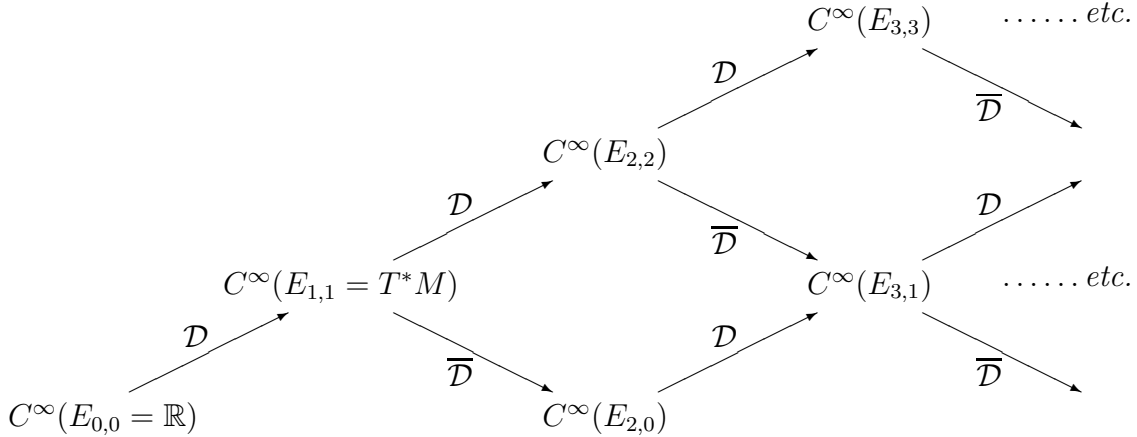
\begin{figure}[h]
\caption{The Quaternionic Double Complex}
\begin{center}
\begin{picture}(400,160)

\put(0,0){\makebox(20,0)[c]{$C^\infty(E_{0,0}=\R)$}}
\put(20,10){\vector(2,1){60}}
\put(45,27){\makebox(0,0)[b]{$\D'$}}

\put(100,50){\makebox(0,0)[c]{$C^\infty(E_{1,1}=T^*M)$}}
\put(120,60){\vector(2,1){60}}
\put(120,40){\vector(2,-1){60}}
\put(145,77){\makebox(0,0)[b]{$\D'$}}
\put(145,23){\makebox(0,0)[t]{$\DD$}}

\put(200,0){\makebox(0,0)[c]{$C^\infty(E_{2,0})$}}
\put(220,10){\vector(2,1){60}}
\put(245,27){\makebox(0,0)[b]{$\D'$}}

\put(200,100){\makebox(0,0)[c]{$C^\infty(E_{2,2})$}}
\put(220,110){\vector(2,1){60}}
\put(220,90){\vector(2,-1){60}}
\put(245,127){\makebox(0,0)[b]{$\D'$}}
\put(245,73){\makebox(0,0)[t]{$\DD$}}

\put(300,50){\makebox(0,0)[c]{$C^\infty(E_{3,1})$}}
\put(320,60){\vector(2,1){60}}
\put(320,40){\vector(2,-1){60}}
\put(345,77){\makebox(0,0)[b]{$\D'$}}
\put(345,23){\makebox(0,0)[t]{$\DD$}}
\put(350,50){\ldots\ldots{\it etc.}}

\put(300,150){\makebox(0,0)[c]{$C^\infty(E_{3,3})$}}
\put(320,140){\vector(2,-1){60}}
\put(345,123){\makebox(0,0)[t]{$\DD$}}
\put(350,150){\ldots\ldots{\it etc.}}

\end{picture}
\end{center}
\end{figure}

Here is our quaternionic analogue of the Dolbeault complex. 
There are strong similarities between this and the real Dolbeault complex (Figure 3.1). 
Again, instead of a diamond as in the Dolbeault complex, the quaternionic version
only extends upwards to form an isosceles triangle. 
This is essentially because there is one irreducible
U(1)-representation 
for each integer, whereas there is one irreducible 
Sp(1)-representation only for each nonnegative integer.

By definition, the bundle $E_{k,k}$ is the bundle $A^k$ of 
(\ref{Akdfn}) --- they are both the subbundle of $\Lambda^kT^*M$ which includes all Sp(1)-representations of the form $V_k$. Thus the leading edge of the double complex  
\[0\lra C^\infty(E_{0,0})\stackrel{\D'}{\lra}
       C^\infty(E_{1,1})\stackrel{\D'}{\lra}
       C^\infty(E_{2,2})\stackrel{\D'}{\lra}\ldots
\stackrel{\D'}{\lra}
       C^\infty(E_{2n,2n})\stackrel{\D'}{\lra}0
\]
is precisely the complex (\ref{Salcomplex}) discovered by Salamon.

\pagebreak
\begin{example}{\bf Four Dimensions}

This double complex is already very well-known and understood in four 
dimensions. Here there is a splitting only in the middle dimension, 
$\Lambda^2T^*M\cong V_2\op 3V_0$. 
Let $I$, $J$ and $K$ be local almost complex structures at $x\in M$, and let $e^0\in T^*_xM$. Let $e^1=I(e^0)$, $e^2=J(e^0)$ and 
$e^3=K(e^0)$. In this way we obtain a basis $\{e^0,\ldots,e^3\}$ for $T^*_xM\cong\H$. Using the notation 
$e^{ijk\ldots}=e^i\wedge e^j\wedge e^k\wedge\ldots${\it etc.}, 
define the 2-forms
\be\omega_1^\pm=e^{01}\pm e^{23},\quad
   \omega_2^\pm=e^{02}\pm e^{31},\quad
   \omega_3^\pm=e^{03}\pm e^{12}.
\ee
Then $I$, $J$ and $K$ all act trivially 
on the $\omega_j^-$, so
$E_{2,0}=\langle\omega_1^-,\omega_2^-,\omega_3^-\rangle$. 
The action of $\sp(1)$ on the $\omega_j^+$ is given by the multiplication table

\be \begin{array}{ccc}
\begin{array}{rcl}
I(\omega_1^+) & = & 0  \\       
J(\omega_1^+) & = & -2\omega_3^+ \\
K(\omega_1^+) & = & 2\omega_2^+ \\
\end{array} 
&
\begin{array}{rcl}
I(\omega_2^+) & = & 2\omega_3^+ \\
J(\omega_2^+) & = & 0 \\
K(\omega_2^+) & = & -2\omega_1^+ 
\end{array}
&
\begin{array}{rcl}
I(\omega_3^+) & = & -2\omega_3^+\\
J(\omega_3^+) & = & 2\omega_1^+ \\
K(\omega_3^+) & = & 0.
\end{array}
\end{array}
\label{sdtable}
\ee    
These are the relations of the irreducible $\sp(1)$-representation 
$V_2$, and we see that 
$E_{2,2}=\langle\omega_1^+,\omega_2^+,\omega_3^+\rangle$.

These bundles will be familiar to most readers: 
$E_{2,2}$ is the bundle of {\it self-dual} 2-forms $\Lambda^2_+$ and $E_{2,0}$ is the bundle of {\it anti-self-dual} 2-forms 
$\Lambda^2_-$.
The celebrated splitting $\Lambda^2T^*M\cong\Lambda^2_+\op\Lambda^2_-$ 
is an invariant of the conformal class 
of any Riemannian 4-manifold, and $I^2+J^2+K^2=-4(\ast+1)$, where 
$\ast:\Lambda^kT^*M\ra\Lambda^{4-k}T^*M$ is the Hodge star map. 

This also serves to explain the special definition that a 4-manifold is said to be quaternionic if it is self-dual and conformal.
The relationship between quaternionic, almost complex and Riemannian 
structures in four dimensions is described in detail in
\cite[Chapter 7]{SalHol}. 

\label{4dimexample}
\end{example}


Because there is no suitable quaternionic version of holomorphic coordinates, there is no `nice' co-ordinate expression for a typical section of $C^\infty(E_{k,r})$. In order to determine the decomposition of a differential form, the simplest way the author
has found is to use the Casimir operator ${\cal C}=I^2+J^2+K^2$. 
Consider a $k$-form $\alpha$. Then $\alpha\in E_{k,r}$ if and only if $(I^2+J^2+K^2)(\alpha)=-r(r+2)\alpha$. This mechanism also allows us to work out expressions for $\D'$ and $\DD$ acting on $\alpha$.

\begin{lem}
Let $\alpha\in C^\infty(E_{k,r})$. Then 
\[\D'\alpha=-\frac{1}{4}\left((r-1)+
            \frac{1}{r+1}(I^2+J^2+K^2)\right)d\alpha\]
and
\[\DD\alpha=\frac{1}{4}\left((r+3)+
            \frac{1}{r+1}(I^2+J^2+K^2)\right)d\alpha.\]
\label{Dformulas}
\end{lem}

\begin{proof}
We have $d\alpha=\D'\alpha+\DD\alpha$, where 
$\D'\alpha\in E_{k+1,r+1}$ and $\DD\alpha\in E_{k+1,r-1}$.
Applying the Casimir operator gives 
\[(I^2+J^2+K^2)(d\alpha)=-(r+1)(r+3)\D'\alpha -(r+1)(r-1)\DD\alpha.\]
Rearranging these equations gives $\D'\alpha$ and $\DD\alpha$. 
\end{proof}

Note that our decomposition is of real- as well as complex- valued  forms; the operators $\D'$ and $\DD$ map real forms to real forms.

Writing $\D'_{k,r}$ for the particular map $\D':C^\infty(E_{k,r})\ra C^\infty(E_{k+1,r+1})$,
we define the quaternionic cohomology groups
\be H^{k,r}_\D'(M)=\frac{\Ker (\D'_{k,r})}
                               {\Im (\D'_{k-1,r-1})}.
\label{cohomdfn}
\ee



\section{Ellipticity and the Double Complex}
\label{Elipsec}

In this section we shall determine where our double complex is elliptic and where it is not. Its properties are extremely like those of the real Dolbeault complex studied earlier: the quaternionic double complex is elliptic everywhere except on the bottom row. Though this is much more difficult to prove for the quaternionic double complex, the fundamental reason is the same as for the real Dolbeault complex: 
it is the isosceles triangle as opposed to diamond shape which causes ellipticity to fail for the bottom row, because $d=\D'$ on $E_{2k,0}$ and the projection from $d(C^\infty(E_{2k,0}))$ to 
$C^\infty(E_{2k+1,1})$ is the identity.

Here is the main result of this section:

\begin{thm}
For $2k\ge 4$, the complex 
\[0\ra
E_{2k,0}\stackrel{\D'}{\ra}
E_{2k+1,1}\stackrel{\D'}{\ra}
E_{2k+2,2}\stackrel{\D'}{\ra}
\ldots\stackrel{\D'}{\ra}
E_{2n+k,2n-k}\stackrel{\D'}{\ra}0
\] 
is elliptic everywhere except at $E_{2k,0}$ and $E_{2k+1,1}$, where it is not elliptic. 

For $k=1$ the complex is elliptic everywhere except at $E_{3,1}$, where it is not elliptic.

For $k=0$ the complex is elliptic everywhere.
\label{elipthm}
\end{thm}

The rest of this section provides a proof of this theorem. 

On a complex manifold $M^{2n}$ with holomorphic coordinates $z^j$,  
the exterior forms
$dz^{a_1}\wedge\ldots\wedge dz^{a_p}\wedge 
d\bar{z}^{b_1}\wedge\ldots\wedge d\bar{z}^{b_q}$
span $\Lambda^{p,q}$. 
This allows us to decompose any form $\omega\in\Lambda^{p,q}$, making it much easier to write down the kernels and images of maps which involve exterior multiplication. 
On a quaternionic manifold $M^{4n}$ there is unfortunately no easy way to write down a local frame for the bundle $E_{k,r}$, because there is no quaternionic version of `holomorphic coordinates'. However, we can decompose $E_{k,r}$ just enough to enable us to prove Theorem \ref{elipthm}.

A principal observation is that since ellipticity is a local property, we can work on $\H^n$ without loss of generality. Secondly, since $\,\GL(n,\H)$\, acts transitively on
 $\,\H^n\setminus\{0\}$, if the symbol 
sequence \,
$\ldots\stackrel{\sigma_{e^0}}{\lra}
E_{k,r}\stackrel{\sigma_{e^0}}{\lra}
E_{k+1,r+1}\stackrel{\sigma_{e^0}}{\lra}\ldots$\, is exact for any 
nonzero $e^0\in T^*\H^n$ then it is exact for all nonzero 
$\xi\in T^*\H^n$.
To prove Theorem \ref{elipthm}, we choose one such $e^0$ and analyse the spaces $E_{k,r}$ accordingly. 


\subsection{Decomposition of the Spaces $E_{k,r}$}
\label{Lieinsec}

Let $e^0\in T^*_x\H^n\cong\H^n$ and let $(I, J, K)$ be the standard hypercomplex structure on $\H^n$. As in Example \ref{4dimexample},
 define $e^1=I(e^0)$, $e^2=J(e^0)$ and 
$e^3=K(e^0)$, so that $\langle e^0,\ldots,e^3\rangle\cong\H$. 
In this way we single out a particular copy of \H\ which we call 
$\H_0$, obtaining a (nonnatural) splitting 
$T^*_x\H^n\cong \H^{n-1}\op\H_0$  which is preserved by action of the hypercomplex structure. 
This induces the decomposition 
$\Lambda^k\H^n\cong\bigoplus_{l=0}^4
                   \Lambda^{k-l}\H^{n-1}\ot\Lambda^l\H_0,$                  
which decomposes each  
$E_{k,r}\subset\Lambda^k\H^n$ according to how many differentials in the $\H_0$-direction are present. 

\begin{dfn}
Define the space $E_{k,r}^l$ to be the subspace of $E_{k,r}$ consisting of exterior forms with precisely $l$ differentials in the 
$\H_0$-direction, {\it i.e.}
\[E_{k,r}^l\equiv E_{k,r}\cap(\Lambda^{k-l}\H^{n-1}\ot\Lambda^l\H_0).\]
Then $E_{k,r}^l$ is preserved by the induced action of the hypercomplex structure on $\Lambda^k\H^n$. Thus we obtain an invariant  decomposition
$E_{k,r}=E_{k,r}^0\op E_{k,r}^1\op E_{k,r}^2\op E_{k,r}^3\op
E_{k,r}^4$. 
Note that we can identify $E_{k,r}^0$ on $\H^n$ with  
$E_{k,r}$ on $\H^{n-1}$.
\label{aspacedfn} 
\end{dfn}

(Throughout the rest of this section, juxtaposition of exterior forms will denote exterior multiplication, for example 
$\alpha e^{ij}$ means $\alpha\wedge e^{ij}$.)

We can decompose these summands still further. Consider, for example, the bundle $E_{k,r}^1$. An exterior form $\alpha\in E_{k,r}^1$ is of the form
$\alpha_0e^0+\alpha_1e^1+\alpha_2e^2+\alpha_3e^3$, where 
$\alpha_j\in\Lambda^{k-1}\H^{n-1}$. 
Thus $\alpha$ is an element of $\Lambda^{k-1}\H^{n-1}\ot 2V_1$, 
since $\H_0\cong 2V_1$ as an $\sp(1)$-representation. Since $\alpha$ is in a copy of the representation $V_r$, it follows from the isomorphism  
$\,V_r\ot V_1\cong V_{r+1}\op V_{r-1}$\, that the $\alpha_j$ must be in a combination of 
$V_{r+1}$ and $V_{r-1}$ representations, {\it i.e.}
$\alpha_j\in E_{k-1,r+1}^0\op E_{k-1,r-1}^0$. We write 
\[\alpha=\alpha^++\alpha^-=
        (\alpha^+_0+\alpha^-_0)e^0+
        (\alpha^+_1+\alpha^-_1)e^1+
        (\alpha^+_2+\alpha^-_2)e^2+
        (\alpha^+_3+\alpha^-_3)e^3,
\]
where $\alpha_j^+\in E_{k-1,r+1}^0$ and $\alpha_j^-\in E_{k-1,r-1}^0$.

The following Lemma allows us to consider $\alpha^+$ and 
$\alpha^-$ separately.

\begin{lem}
If $\alpha=\alpha^++\alpha^-\in E_{k,r}^1$ then both 
$\alpha^+$ and $\alpha^-$ must be in $E_{k,r}^1$.

\label{Lieinlem1}
\end{lem}

\begin{proof}
In terms of representations, the situation is of the form 
\[ (pV_{r+1}\op qV_{r-1})\ot 2V_1\cong 
              2p(V_{r+2}\op V_r)\op 2q(V_r\op V_{r-2}),
\]
where $\alpha^+\in pV_{r+1}$ and $\alpha^-\in qV_{r-1}$.
For $\alpha$ to be in the representation $2(p+q)V_r$, its components in the representations $2pV_{r+2}$ and $2qV_{r-2}$ must both vanish separately. The component in $2pV_{r+2}$ comes entirely from 
$\alpha^+$, so for this to vanish we must have $\,\alpha^+\in 2(p+q)V_r$\, independently of $\alpha^-$. 
Similarly, for the component in $2qV_{r-2}$ to vanish, we must have 
$\,\alpha^-\in 2(p+q)V_r$.
\end{proof}

Thus we decompose the space $E_{k,r}^1$ into two summands, one coming from $E_{k-1,r-1}^0\ot 2V_1$ and one from $E_{k-1,r+1}^0\ot 2V_1$.
We extend this decomposition to the cases $l=0,2,3,4$, defining the following notation.

\begin{dfn}
Define the bundle $E_{k,r}^{l,m}$ to be the subbundle of $E_{k,r}^l$ arising from $V_m$-type representations in $\Lambda^{k-l}\H^{n-1}$.
In other words, 
\[E_{k,r}^{l,m}\equiv (E_{k-l,m}^0\ot\Lambda^l\H_0)\cap E_{k,r}.\]
\label{abspacedfn}
\end{dfn}

To recapitulate: for the space $E_{k,r}^{l,m}$, the bottom left index $\,k$\, refers to the exterior power of the form 
$\alpha\in\Lambda^k\H^n$; the bottom right index $\,r$\, refers to the irreducible Sp(1)-representation in which $\alpha$ lies; the top left index $\,l$\, refers to the number of differentials in the 
$\H_0$-direction and the top right index $\,m$\, refers to the irreducible Sp(1)-representation of the contributions 
from $\Lambda^{k-a}\H^{n-1}$\, {\it before} wedging with forms in the 
$\H_0$-direction.
This may appear slightly fiddly: it becomes rather simpler when we consider the specific splittings which Definition \ref{abspacedfn}
allows us to write down.

\begin{lem}
Let $E_{k,r}^{l,m}$ be as above. We have the following decompositions:
\[E_{k,r}^0=E_{k,r}^{0,r}\quad\quad 
  E_{k,r}^1=E_{k,r}^{1,r+1}\op E_{k,r}^{1,r-1}\quad\quad
  E_{k,r}^2=E_{k,r}^{2,r+2}\op E_{k,r}^{2,r}\op E_{k,r}^{2,r-2}
\]
\[E_{k,r}^3=E_{k,r}^{3,r+1}\op E_{k,r}^{3,r-1}\quad\quad\mbox{and}
                                              \quad\quad
  E_{k,r}^4=E_{k,r}^{4,r}.
\]
\label{Lieinlem2}
\end{lem}

\begin{proof}
The first isomorphism is trivial, as is the last (since the hypercomplex structure acts trivially on 
$\Lambda^4\H_0$).
The second isomorphism is Lemma \ref{Lieinlem1}, and the fourth follows in exactly the same way since $\Lambda^3\H_0\cong 2V_1$ also. The middle isomorphism follows a similar argument.
\end{proof}

Recall the self-dual forms and anti-self-dual forms in 
Example \ref{4dimexample}. 
The bundle $E^{2,r}_{k,r}$ splits according to whether its contribution from $\Lambda^2\H_0$ is self-dual or anti-self-dual. We will call these summands $E_{k,r}^{2,r+}$ and 
$E_{k,r}^{2,r-}$ respectively, so 
$E^{2,r}_{k,r}=E_{k,r}^{2,r+}\op E_{k,r}^{2,r-}$. 

\subsection{Lie in conditions}

We have analysed the bundle $E_{k,r}$ into a number of different subbundles. We now determine when a
particular exterior form lies in one of these subbundles. 
Consider a form 
$\alpha=\alpha_1e^{s_1\ldots s_a}+\alpha_2e^{t_1\ldots t_a}+\ldots
$ {\it etc.} where $\alpha_j\in E_{k-a,b}^0$.
For $\alpha$ to lie in one of the spaces $\,E_{k,r}^{a,b}$\, the 
$\alpha_j$ will usually have to satisfy some simultaneous equations.
Since these are the conditions for a form to lie in a particular Lie algebra representation, we will refer to such equations as `Lie In Conditions'. 

To begin with, we mention three trivial lie in conditions. Let 
$\alpha\in E_{k,r}^0$. That $\alpha\in E_{k,r}^{0,r}$ is obvious, as is \,$\alpha e^{0123}\in E_{k,r}^{4,r}$, since wedging with 
$e^{0123}$ has no effect on the $\sp(1)$-action. 
Likewise, the $\sp(1)$-action on the anti-self-dual 2-forms 
\,$\omega_1^-=e^{01}-e^{23}$, $\omega_2^-=e^{02}-e^{31}$ and
$\omega_3^-=e^{03}-e^{12}$ is trivial, so 
$\alpha\omega_j^-\in E_{k,r}^{2,r-}$ for all $j=1,2,3.$ 

This leaves the following three situations: those arising from taking exterior products with 1-forms, 3-forms and the self-dual 2-forms 
$\omega_j^+$.  
As usual when we want to know which representation an exterior form is in, we apply the Casimir operator.

\subsubsection{The cases $l=1$ and $l=3$}

Let \,$\alpha_j\in E_{k,r}^0$. Then 
\,$\alpha=\alpha_0e^0+\alpha_1e^1+\alpha_2e^2+\alpha_3e^3\in 
E_{k+1,r+1}^{1,r}\op E_{k+1,r-1}^{1,r}$, and 
$\alpha$ is entirely in $E_{k+1,r+1}^{1,r}$ if and only if
$\,(I^2+J^2+K^2)\alpha=-(r+1)(r+3)\alpha$.

By the usual (Leibniz) rule for a Lie algebra action on a tensor product, we have that 
\linebreak
$I^2(\alpha_je^j)=I^2(\alpha_j)e^j+2I(\alpha_j)I(e^j)+
                      \alpha_j I^2(e^j)$, {\it etc}.
Thus
\begin{eqnarray} (I^2+J^2+K^2)\alpha
& = & \sum_{j=0}^3\Big[(I^2+J^2+K^2)(\alpha_j)e^j 
       + \alpha_j(I^2+J^2+K^2)(e^j)\ +            \nonumber\\
&   & \quad\quad\quad\quad +\ 2\Big(I(\alpha_j)I(e^j) + 
      J(\alpha_j)J(e^j) + K(\alpha_j)K(e^j)\Big) 
      \Big] \nonumber \\
& = & -r(r+2)\alpha - 3\alpha + 2\sum_{j=0}^3 
       \Big(I(\alpha_j)I(e^j) +
        J(\alpha_j)J(e^j) +
        K(\alpha_j)K(e^j)\Big) \nonumber\\
& = & (-r^2-2r-3)\alpha + 
2\begin{array}[t]{l}
   \Big(I(\alpha_0)e^1-I(\alpha_1)e^0+I(\alpha_2)e^3-I(\alpha_3)e^2+  
\\  
   +J(\alpha_0)e^2-J(\alpha_1)e^3-J(\alpha_2)e^0+J(\alpha_3)e^1+ 
\\
   +K(\alpha_0)e^3+K(\alpha_1)e^2-K(\alpha_2)e^1-K(\alpha_3)e^0\Big).
\end{array} 
\nonumber \\
\label{Casalphaeqn1}
\end{eqnarray}
For $\alpha\in E_{k+1,r+1}^{1,r}$ we need this to be equal to 
$-(r+1)(r+3)\alpha$, which is the case if and only if 
\begin{eqnarray*} -r\alpha & = &
   I(\alpha_0)e^1-I(\alpha_1)e^0+I(\alpha_2)e^3-I(\alpha_3)e^2 
  +J(\alpha_0)e^2-J(\alpha_1)e^3-J(\alpha_2)e^0+J(\alpha_3)e^1+ 
\nonumber \\
&  & +\ K(\alpha_0)e^3+K(\alpha_1)e^2-K(\alpha_2)e^1-K(\alpha_3)e^0.
\end{eqnarray*}

Since the $\alpha_j$ have no $e^j$-factors and the action of 
$I$, $J$ and $K$ preserves this property, this equation can only be satisfied if 
it holds for each of the $e^j$-components separately. It follows that $\alpha\in E_{k+1,r+1}^{1,r}$ if and only if $\alpha_0$, $\alpha_1$,
$\alpha_2$ and $\alpha_3$ satisfy the following lie in conditions:
\footnote{Our interest in these conditions arises from a consideration 
of exterior forms, but the equations describe $\sp(1)$-representations in general: they are the conditions that 
$\alpha\in V_r\ot V_1$ must satisfy to be in the $V_{r+1}$ subspace of 
$V_{r+1}\op V_{r-1}\cong V_r\ot V_1$. The other lie in conditions have similar interpretations.} 
 
\be \begin{array}{rcl}
r\alpha_0 - I(\alpha_1) - J(\alpha_2) - K(\alpha_3) & = & 0 \\
r\alpha_1 + I(\alpha_0) + J(\alpha_3) - K(\alpha_2) & = & 0 \\
r\alpha_2 - I(\alpha_3) + J(\alpha_0) + K(\alpha_1) & = & 0 \\
r\alpha_3 + I(\alpha_2) - J(\alpha_1) + K(\alpha_0) & = & 0. \\ 
\end{array}
\label{LICeqn1}
\ee

Suppose instead that $\alpha\in E_{k+1,r-1}^{1,r}$. Then 
$(I^2+J^2+K^2)\alpha=-(r-1)(r+1)\alpha$. Putting this alternative into 
Equation (\ref{Casalphaeqn1}) gives the result that 
$\alpha\in E_{k+1,r-1}^{1,r}$ if and only if 

\be \begin{array}{rcl}
(r+2)\alpha_0 + I(\alpha_1) + J(\alpha_2) + K(\alpha_3) & = & 0 \\
(r+2)\alpha_1 - I(\alpha_0) - J(\alpha_3) + K(\alpha_2) & = & 0 \\
(r+2)\alpha_2 + I(\alpha_3) - J(\alpha_0) - K(\alpha_1) & = & 0 \\
(r+2)\alpha_3 - I(\alpha_2) + J(\alpha_1) - K(\alpha_0) & = & 0. \\ 
\end{array}
\label{LICeqn2}
\ee

Consider now 
$\alpha=\alpha_0e^{123}+\alpha_1e^{032}+\alpha_2e^{013}+
\alpha_3e^{021}\in E_{k+3,r+1}^{3,r}\op E_{k+3,r-1}^{3,r}$. 
Since $\Lambda^3\H_0\cong \H_0$, the lie in conditions are exactly the same: for 
$\alpha$ to be in $E_{k+3,r+1}^{3,r}$ we need the $\alpha_j$ to satisfy Equations (\ref{LICeqn1}), and for 
$\alpha$ to be in $E_{k+3,r-1}^{3,r}$ we need the $\alpha_j$ to satisfy Equations (\ref{LICeqn2}).


\subsubsection{The case $l=2$}

We have already noted that wedging a form $\beta\in E_{k,r}^0$ with an anti-self-dual 2-form $\omega_j^-$ has no effect on the $\sp(1)$-action, so $\beta\omega_j^-\in E_{k+2,r}^{2,r-}$. Thus we only have to consider the effect of wedging with the self-dual 2-forms 
$\,\langle\omega_1^+,\omega_2^+,\omega_3^+\rangle
\cong V_2\subset\Lambda^2\H_0$. By the Clebsch-Gordon formula, the decomposition takes the form 
$\,V_r\ot V_2\cong V_{r+2}\op V_r\op V_{r-2}$. Thus for 
$\beta=\beta_1\omega_1^+ +\beta_2\omega_2^+ +\beta_3\omega_3^+$ we want to establish the lie in conditions for $\beta$ to be in 
$E_{k+2,r+2}^{2,r}$, $E_{k+2,r}^{2,r+}$ and $E_{k+2,r-2}^{2,r}$.  

We calculate these lie in conditions in a similar fashion to the previous cases, by considering the action of the Casimir operator 
$I^2+J^2+K^2$ on $\beta$ and using the multiplication table (\ref{sdtable}). The following lie in conditions are then easy to deduce:

\be
\beta\in E_{k+2,r+2}^{2,r}\Longleftrightarrow
\left\{   \begin{array}{rcl}
             (r+4)\beta_1 & = & J(\beta_3)-K(\beta_2)\\
             (r+4)\beta_2 & = & K(\beta_1)-I(\beta_3)\\
             (r+4)\beta_3 & = & I(\beta_2)-J(\beta_1).\\
\end{array}\right.
\label{LIC1}
\ee

\be
\beta\in E_{k+2,r}^{2,r+}\Longleftrightarrow
\left\{   \begin{array}{rcl}
             2\beta_1 & = & J(\beta_3)-K(\beta_2)\\
             2\beta_2 & = & K(\beta_1)-I(\beta_3)\\
             2\beta_3 & = & I(\beta_2)-J(\beta_1).\\
\end{array}\right.
\label{LIC2}
\ee

\be
\beta\in E_{k+2,r-2}^{2,r}\Longleftrightarrow
\left\{   \begin{array}{rcl}
             (2-r)\beta_1 & = & J(\beta_3)-K(\beta_2)\\
             (2-r)\beta_2 & = & K(\beta_1)-I(\beta_3)\\
             (2-r)\beta_3 & = & I(\beta_2)-J(\beta_1).\\
\end{array}\right.
\label{LIC3}
\ee

Equation (\ref{LIC2}) is particularly interesting. Since this 
equation singles out the 
\linebreak
$V_r$-representation in the direct sum 
$\,V_{r+2}\op V_r\op V_{r+2}$, it must have $\,\dim V_r=r+1$\, linearly independent solutions.
Let $\beta_0\in V_r$ and let $\beta_1=I(\beta_0)$,
$\beta_2=J(\beta_0)$, $\beta_3=K(\beta_0)$.
Using the Lie algebra relations $2I=[J,K]=JK-KJ$, it is easy to see that $\beta_1$, $\beta_2$ and $\beta_3$ satisfy Equation \ref{LIC2}.
Moreover, there are $\,r+1$\, linearly independent solutions of this form                                 
 (for $r\neq 0$).  We conclude that {\it all} the solutions of Equation (\ref{LIC2}) take the form 
$\beta_1=I(\beta_0)$,
$\beta_2=J(\beta_0)$, $\beta_3=K(\beta_0)$.


\subsection{The Symbol Sequence and Proof of Theorem \ref{elipthm}}

We now describe the principal symbol of $\D'$, and examine its behaviour in the context of the decompositions of Definition 
\ref{aspacedfn} and Lemma \ref{Lieinlem2}. This leads to a proof of Theorem  \ref{elipthm}. 
First we obtain the principal symbol from the formula for $\D'$ in Lemma \ref{Dformulas} by replacing 
$d\alpha$ with $\alpha e^0$. 

\begin{prop}
Let $\,x\in \H^n$, $e^0\in T^*_x\H^n$\, and $\,\alpha\in E_{k,r}$. The principal symbol mapping \,
$\sigma_{\D'}(x,e^0):E_{k,r}\ra E_{k+1,r+1}$ 
is given by 
\[\sigma_{\D'}(x,e^0)(\alpha)=\frac{1}{2(r+1)}\left(
        (r+2)\alpha e^0-I(\alpha)e^1-J(\alpha)e^2-K(\alpha)e^3
                     \right).
\]
\label{symbolformula}
\end{prop}

\begin{proof}
Replacing $d\alpha$ with $\alpha e^0$ in the formula for $\D'$ obtained in Lemma \ref{Dformulas}, we have
\begin{eqnarray*}\sigma_\D'(x,e^0)(\alpha) & = & -\frac{1}{4}\Big(
               (r-1)+\frac{1}{r+1}(I^2+J^2+K^2)\Big)\alpha e^0 \\
& = & \frac{-1}{4(r+1)}\Big[\left(\,(r-1)(r+1)-r(r+2)-3\,\right)
                             \alpha e^0 \\ 
&   & \ \ \ \ \ \ \ \ \ \ \ \ 
       + 2\left(I(\alpha)e^1+J(\alpha)e^2+K(\alpha)e^3\right)\Big]\\
& = & \frac{1}{2(r+1)}\Big((r+2)\alpha e^0-I(\alpha)e^1-J(\alpha)      
               e^2-K(\alpha)e^3\Big),
\end{eqnarray*}
as required.
\end{proof}

\begin{cor}
The principal symbol $\,\sigma_\D'(x,e^0)$\, maps the space 
$\,E_{k,r}^{l,m}$\, to the space $\,E_{k+1,r+1}^{l+1,m}$.
\end{cor}

\begin{proof}
We already know that $\sigma_\D':E_{k,r}\ra E_{k+1,r+1}$, by definition. 
Using Lemma \ref{symbolformula}, we see that $\sigma_\D'(x,e^0)$ increases the number of differentials in the $\H_0$-direction by one, so the index $l$ increases by one.  The only action in the other directions is the $\sp(1)$-action, which preserves the irreducible decomposition of the contribution from $\Lambda^{k-a}\H^{n-1}$, so the index $m$ remains the same. 
\end{proof}

To save space we shall use $\sigma$ as an abbreviation for 
$\sigma_\D'(x,e^0)$ for the rest of this section. 
The point of all this work on decomposition now becomes apparent. Since $\sigma:E_{k,r}^l\ra E_{k+1,r+1}^{l+1}$, we can reduce the indefinitely long symbol sequence 
\[\ldots
\stackrel{\sigma}{\lra} E_{k-1,r-1} 
\stackrel{\sigma}{\lra} E_{k,r} 
\stackrel{\sigma}{\lra} E_{k+1,r+1} 
\stackrel{\sigma}{\lra} \ldots 
\mbox{{\it etc.}}
\]
to the 5-space sequence 
\be \ \ \ \ \ \ \ \ \ 0\stackrel{\sigma}{\lra} E_{k-2,r-2}^0\stackrel{\sigma}{\lra}
E_{k-1,r-1}^1\stackrel{\sigma}{\lra} 
E_{k,r}^2\stackrel{\sigma}{\lra} 
E_{k+1,r+1}^3\stackrel{\sigma}{\lra} 
E_{k+2,r+2}^4\stackrel{\sigma}{\lra} 0.
\label{5spacesequence}
\ee
Using Lemma \ref{Lieinlem2} as well, we can analyse this sequence still further according to the different (top right) $m$-indices, obtaining three short sequences (for $k\ge 2$, $k\equiv r \bmod 2$)
\be
\begin{array}{ccccccccccccc}
& & & & 0 & \ra & E_{k,r}^{2,r+2}& \ra & E_{k+1,r+1}^{3,r+2}& 
            \ra & E_{k+2,r+2}^{4,r+2} & \ra & 0 \\
& & & &   &      & \op && \op \\
& & 0 &\ra & E_{k-1,r-1}^{1,r} & \ra & E_{k,r}^{2,r} & \ra
                 & E_{k+1,r+1}^{3,r} & \ra & 0 \\
& &   &     & \op & & \op \\
0 &\ra & E_{k-2,r-2}^{0,r-2} & \ra & E_{k-1,r-1}^{1,r-2} & \ra
                 & E_{k,r}^{2,r-2} & \ra & 0.\\
\end{array}
\label{splitarray}
\ee
This reduces the problem of determining where the operator  $\D'$ is elliptic to the problem of determining when these three sequences are exact.

For a sequence $0\ra A\ra B\ra C\ra 0$ to be exact, it is necessary that $\dim A-\dim B+\dim C=0$. Given this condition, if the sequence is exact at any two out of $A$, $B$ and $C$ it is exact at the third.
We shall show that for $r\neq 0$ this dimension sum does equal zero.

\begin{lem}
For $r>0$, each of the sequences in (\ref{splitarray}) satisfies the dimension condition above, {\it i.e.} the alternating sum of the dimensions vanishes.
\label{dimsumlem}
\end{lem}

\begin{proof} Let $r>0$.
We calculate the dimensions of the spaces $E_{k,r}^{l,m}$ for 
$l=0,\ldots,4$. Recall the notation $E_{k,r}=\epsilon^n_{k,r}V_r$ from Definition \ref{Enotationdfn}. It is clear that 
$\,\dim E_{k,r}^0=(r+1)\epsilon^{n-1}_{k,r}$, since 
$E_{k,r}^0$ on $\H^n$ is sinply $E_{k,r}$ on $\H^{n-1}$. 
Thus $\,\dim E_{k-2,r-2}^{0,r-2}=(r-1)\epsilon^{n-1}_{k-2,r-2}$\, and
$\,\dim E_{k+2,r+2}^{4,r+2}=(r+3)\epsilon^{n-1}_{k-2,r+2}$. 

The cases $a=1$ and $a=3$ are easy to work out since they are of the form $E_{k,r}^0\ot 2V_1$.
For $a=1$, we have 
$\,\dim E_{k-1,r-1}^{1,r-2}=2r\epsilon^{n-1}_{k-2,r-2}$\, and 
$\,\dim E_{k-1,r-1}^{1,r}=2r\epsilon^{n-1}_{k-2,r}$.
For $a=3$, $\dim E_{k+1,r+1}^{3,r}=2(r+2)\epsilon^{n-1}_{k-2,r}$ and
$\,\dim E_{k+1,r+1}^{3,r+2}=2(r+2)\epsilon^{n-1}_{k-2,r+2}$. 

The case $a=2$ is slightly more complicated, as we have to take into account exterior products with the self-dual 2-forms $V_2$ and anti-self-dual 2-forms $V_0$ in  
$\Lambda^2\H_0$. 
The spaces $E^{2,r+2}_{k,r}$ and $E^{2,r-2}_{k,r}$ receive contributions only from the self-dual part $V_2$, from which we infer that $\,\dim E_{k,r}^{2,r+2}=(r+1)\epsilon^{n-1}_{k-2,r+2}$\, and 
$\,\dim E_{k,r}^{2,r-2}=(r+1)\epsilon^{n-1}_{k-2,r-2}$. 
Finally, the space $E_{k,r}^{2,r+}$ has dimension 
$\,(r+1)\epsilon^{n-1}_{k-2,r}$\, 
and the space $E_{k,r}^{2,r-}$ has dimension $\,3(r+1)\epsilon^{n-1}_{k-2,r}$, giving $E_{k,r}^{2,r}$ a total dimension of 
$\,4(r+1)\epsilon^{n-1}_{k-2,r}$.

It is now a simple matter to verify that for the top sequence of (\ref{splitarray}) 
\[\epsilon^{n-1}_{k-2,r+2}(r+1-2(r+2)+r+3)=0,\]
for the middle sequence
\[\epsilon^{n-1}_{k-2,r}(2r-4(r+1)+2(r+2))=0,\]
and for the bottom sequence
\[\epsilon^{n-1}_{k-2,r-2}(r-1-2r+r+1)=0.\]
\end{proof}

The case $r=0$ is different. Here the bottom sequence of 
(\ref{splitarray}) disappears 
altogether, the top sequence still being exact. 
Exactness is lost in the middle sequence. Since the isomorphism 
$\epsilon^{n-1}_{k-2,0}V_0\ot V_2\cong \epsilon^{n-1}_{k-2,0}V_2$\, gives no trivial $V_0$-representations, there is no space $E_{k,0}^{2,0+}$. Thus $E_{k,0}^{2,0}$ is `too small' --- we are left with a sequence 

\[0\lra 3\epsilon^{n-1}_{k-2,0}V_0\lra 
        2\epsilon^{n-1}_{k-2,0}V_1\lra 0,\]
which cannot be exact.
(As there is no space $E_{0,0}^2$, 
this problem does not arise for the leading edge 
\,$0\ra E_{0,0}\ra E_{1,1}\ra\ldots$ {\it etc.})

We are finally in a position to prove Theorem \ref{elipthm}, which now follows from:

\begin{prop}
When $r\neq 0$, the three sequences of (\ref{splitarray}) are exact.  
\label{exactnessprop}
\end{prop}

\begin{proof}
Consider first the top sequence
\,$0\stackrel{\sigma}{\lra}E_{k,r}^{2,r+2}
  \stackrel{\sigma}{\lra}E_{k+1,r+1}^{3,r+2}
  \stackrel{\sigma}{\lra}E_{k+2,r+2}^{4,r+2} \lra 0.$ 
The Clebsch-Gordon formula shows that 
there are no spaces $E_{k+1,r-1}^{3,r+2}$ or 
$E_{k+2,r}^{4,r+2}$. Thus $\DD=0$ on $E_{k,r}^{2,r+2}$ and 
$E_{k+1,r+1}^{3,r+2}$, so $\D'=d$ for the top sequence. It is easy to check using the relevant lie in conditions that the map 
$\wedge e^0:E_{k,r}^{2,r+2}\ra E_{k+1,r+1}^{3,r+2}$ is injective and the map  
$\wedge e^0:E_{k+1,r+1}^{3,r+2}\ra E_{k+2,r+2}^{4,r+2}$ is surjective.

To show exactness at $E_{k-1,r-1}^1$, consider 
$\alpha=\alpha_0e^0+\alpha_1e^1+\alpha_2e^2+\alpha_3e^3 \in
             E_{k-1,r-1}^1$. 
A calculation using Proposition \ref{symbolformula} shows that 
\begin{eqnarray}\sigma(\alpha) & = & \frac{1}{2r}\left(
        (r\alpha_1+I(\alpha_0))e^{10}+
        (r\alpha_2+J(\alpha_0))e^{20}+
        (r\alpha_3+K(\alpha_0))e^{30}+\right. 
\label{symbolzeroeqn} \\
& & \hspace{-.3in}
+\left.(2\alpha_1-J(\alpha_3)+K(\alpha_2))e^{32}+
 (2\alpha_2-K(\alpha_1)+I(\alpha_3))e^{13}+
 (2\alpha_3-I(\alpha_2)+K(\alpha_1))e^{21}\right).\nonumber
\end{eqnarray}
Since the $\alpha_i$ have no $e^j$-components, $\,\sigma(\alpha)=0$\, if and only if all these components vanish. This occurs if and only if 
$\,\alpha_1=-\frac{1}{r}I(\alpha_0)$, 
$\,\alpha_2=-\frac{1}{r}J(\alpha_0)$, 
$\,\alpha_3=-\frac{1}{r}K(\alpha_0)$ (since as remarked in Section 
\ref{Lieinsec} these equations also guarantee that 
$\,2\alpha_1-J(\alpha_3)+K(\alpha_2)=0$ {\it etc.}), in which case it is clear that 
\[\sigma(\alpha)=0\Longleftrightarrow 
\alpha=\sigma\left(\textstyle{\frac{2(r+1)}{r}}\alpha_0\right).\]
This shows that the sequence \,$E^0_{k-2,r-2}\ra E^1_{k-1,r-1}\ra 
E^2_{k,r}$\, is exact. Restricting to $E^{1,r}_{k-1,r-1}$ and 
$E^{1,r-2}_{k-1,r-1}$, we see that exactness holds at these spaces in the middle and bottom sequences respectively of (\ref{splitarray}).

Consider $\alpha\in E_{k-2,r-2}^0$. Then 
\[\sigma(\alpha)=\frac{1}{2(r-1)}\left(r\alpha e^0-
         I(\alpha) e^1 - J(\alpha) e^2 - K(\alpha) e^3\right).
\]
Since these are linearly independent, $\sigma(\alpha)=0$ if and only if 
$\alpha=0$, and $\,\sigma:E_{k-2,r-2}^{0,r-2}\ra E_{k-1,r-1}^{1,r-2}$\, 
is injective.
Hence the bottom sequence
\,$0\lra E_{k-2,r-2}^{0,r-2}
  \stackrel{\sigma}{\lra}E_{k-1,r-1}^{1,r-2}
  \stackrel{\sigma}{\lra}E_{k,r}^{2,r-2} \lra 0$\, is exact. 

Finally, we show that the middle sequence 
\,$0\lra E_{k-1,r-1}^{1,r}
  \stackrel{\sigma}{\lra}E_{k,r}^{2,r}
  \stackrel{\sigma}{\lra}E_{k+1,r+1}^{3,r} \lra 0$\, is exact at 
$E_{k,r}^{2,r}$, which is now sufficient to show that the sequence is exact.  

Let $\beta=\beta_1\omega_1^+
          +\beta_2\omega_2^+
          +\beta_3\omega_3^+\in E_{k,r}^{2,r+}$. 
Recall the lie in condition (\ref{LIC2}) that $\beta$ must take the form $\beta=\frac{1}{r}\left(I(\beta_0)\omega_1^+ + 
             J(\beta_0)\omega_2^+ + 
             K(\beta_0)\omega_3^+\right)$ 
for some $\beta_0\in E_{k-2,r}^0$. (The $\frac{1}{r}$-factor makes no difference here and is useful for cancellations.) Thus a general element of  
$E_{k,r}^{2,r}$ is of the form 
\[\beta+\gamma=\frac{1}{r}
       \left( I(\beta_0)\omega_1^+ + J(\beta_0)\omega_2^+ + 
             K(\beta_0)\omega_3^+ \right) 
     + \gamma_1\omega_1^- + \gamma_2\omega_2^- + \gamma_3\omega_3^-,
\] 
for $\beta_0, \gamma_j\in E_{k-2,r}^{0,r}$. A similar calculation to that of (\ref{symbolzeroeqn}) shows that 
\[ \sigma(\beta+\gamma)=0 \Longleftrightarrow
\left\{\begin{array}{rcl} 
        (r+2)\beta_0+I(\gamma_1)+J(\gamma_2)+K(\gamma_3) & = & 0 \\
        (r+2)\gamma_1-I(\beta_0)-J(\gamma_3)+K(\gamma_2) & = & 0 \\
        (r+2)\gamma_2+I(\gamma_3)-J(\beta_0)-K(\gamma_1) & = & 0 \\
        (r+2)\gamma_3-I(\gamma_2)+J(\gamma_1)-K(\beta_0) & = & 0.
\end{array}\right.
\]
But this is exactly the lie in condition (\ref{LICeqn2}) which we need for $\beta_0e^0+\gamma_1e^1+\gamma_2e^2+\gamma_3e^3$ to be in 
      $E_{k-1,r-1}^{1,r},$ in which case we have 
\[\beta+\gamma=
\sigma\left(2(\beta_0e^0+\gamma_1e^1+\gamma_2e^2+\gamma_3e^3)\right).
\]  
This demonstrates exactness at $E_{k,r}^{2,r}$ and so the middle sequence is exact.
\end{proof}

As a counterexample for the case $\,r=0$ and $k\geq 4$, consider 
$\alpha\in E_{k-4,0}^0$. Then $\alpha e^{0123}\in E_{k,0}^{4,0}$ and 
$\sigma(\alpha e^{0123})=0$, so $\,\sigma:E_{k,0}\ra E_{k+1,1}$ is not injective, which is exactly the same as saying that the symbol sequence is not exact at $E_{k,0}$. It is easy to see that this counterexample does not arise when $k=0$ or $2$, and to show that the maps 
$\,\sigma:E_{0,0}\ra E_{1,1}$\, and $\,\sigma:E_{2,0}\ra E_{3,1}$\, {\it are} injective.

As a counterexample for the case $\,r=1$ and $k\geq 2$, consider 
$\alpha\in E_{k-2,0}^0$.
Then $\alpha e^{123}\in E_{k+1,1}^{3,0}$ 
and $\alpha e^{123}\wedge e^0\in E_{k+2,0}^4$. Thus 
$\sigma(\alpha e^{123})=0$. Since $\alpha e^{123}$ has no 
$\,e^0$-components at all it is clear that 
$\alpha e^{123}\neq \sigma(\beta)$\, for any $\,\beta\in E_{k,0}^2$. Thus 
the symbol sequence fails to be exact at $E_{k+1,1}$.
Again, it is easy to see that this counterexample does not arise  when $k=0$, and to show that the sequence 
$\,E_{0,0}\stackrel{\sigma}{\lra}E_{1,1}
        \stackrel{\sigma}{\lra}E_{2,2}$\, 
{\it is} exact at $E_{1,1}$.

This concludes our proof of Theorem \ref{elipthm}.



\section{Quaternion-valued forms on Hypercomplex Manifolds}
\label{Qformsec}

Let $M$ be a hypercomplex manifold. Then $M$ has a triple $(I,J,K)$ of complex structures which generates the $\sp(1)$-action on $\Lambda^kT^*M$ and which we can identify globally with the imaginary quaternions. 
Joyce has used this identification to define `quaternionic  holomorphic functions', which he calls q-holomorphic functions. A quaternion-valued function 
$f=f_0+f_1i+f_2j+f_3k$ is defined to be q-holomorphic if it satisfies    
a quaternionic version of the Cauchy-Riemann equations \cite[3.3]{Jqalg}
\be df_0+I(df_1)+J(df_2)+K(df_3)=0.
\label{Qholdef}
\ee

This equation can also be obtained by comparing the Sp(1)-representations on $T^*M$ and on the quaternions themselves.
Recall that Equation (\ref{Hnrepeqn}) describes the 
$\Sp(1)\GL(n,\H)$-representation on $\H^n$ as $V_1\ot E$. In the case $n=1$ this reduces to the representation 
\be\H\cong V_1\ot V_1,
\label{Hrepeqn}
\ee 
interpreting the left-hand copy of $V_1$ as the left-action
$(p,q)\mapsto pq$, and the right-hand copy of $V_1$ as the right-action $(p,q)\mapsto qp^{-1}$, for $q\in\H$ and $p\in\Sp(1)$.

We can now use our globally defined hypercomplex structure to combine  the Sp(1)-actions on
$\H$ and $\Lambda^kT^*M$. 
Consider, for example, quaternion-valued exterior forms in the bundle 
$E_{k,r}=\epsilon^n_{k,r}V_r$. The Sp(1)-action on these forms is described by the representation 
\[\H\ot E_{k,r}\cong V_1\ot V_1\ot \epsilon^n_{k,r}V_r.
\]
Leaving the left \H-action untouched, we consider the effect of the right \H-action and the hypercomplex structure simultaneously. This amounts to applying the operators 
\[{\cal I}:\alpha\ra I(\alpha)-\alpha i, \quad\quad
  {\cal J}:\alpha\ra J(\alpha)-\alpha j 
\quad\quad \mbox{and} \quad\quad
  {\cal K}:\alpha\ra K(\alpha)-\alpha k
\]
to $\alpha\in\H\ot E_{k,r}$.
Under this diagonal action the tensor product 
$V_1\ot \epsilon^n_{k,r}V_r$ splits, giving the representation 
\be \H\ot E_{k,r}\cong V_1\ot\epsilon^n_{k,r}(V_{r+1}\op V_{r-1}).
\label{HEsplit}
\ee
This gives a quaternion-valued version of the double complex which has certain advantages over its real-valued counterpart --- for example, in 4-dimensions the whole quaternionic double complex is elliptic. 

Joyce's equation (\ref{Qholdef}) turns out to be one example of this: it is the condition necessary for 
$df$ to lie in the $V_2$-summand of the splitting 
\[\H\ot T^*M\cong V_1\ot 2n(V_2\op V_0).\] 

Joyce's paper also develops a theory of quaternionic algebra based upon left \H-modules whose structure is `augmented' by singling out a particular real subspace. The author has shown that the most important class of these `augmented \H-modules' arises from 
Sp(1)-representations using splittings of the form given by Equation (\ref{HEsplit}). 
This point of view turns out to be very fruitful: it both improves our understanding of Joyce's quaternionic algebra and shows how to apply his theory to many naturally occuring vector bundles over hypercomplex manifolds. This leads not only to Joyce's q-holomorphic functions, but also to quaternionic analogues of holomorphic $k$-forms, the holomorphic tangent and cotangent spaces, and even complex Lie groups and Lie algebras. The author hopes to make this type of quaternionic analysis on hypercomplex manifolds the subject of a future paper.


\end{document}